\def \NN {\mathbb N}
\def \CC {\mathbb C}
\def \QQ {\mathbb Q}
\def \RR {\mathbb R}
\def \ZZ {\mathbb Z}

\def \epsilon{\varepsilon}

\def \K  {{\mathcal K}}

\def \I  {{\mathcal I}}

\def \d {\text{d}}
\def \HH {\Bbb{H}}

\def \half {{\textstyle{ 1\over 2}}}
\def \fine {{\hfill \qedsymbol}}

\def \ep {\epsilon}

\def \si {\sigma}

\def \Log {\text{Log}}

\newcommand{\res}{\text{res}}
\renewcommand{\S}{{\mathcal S}}

\documentclass[12pt,reqno]{amsart}

\usepackage{amsmath,amssymb}

\usepackage{url} 

\numberwithin{equation}{section}

\begin{document}


\title[]{The standard twist of $L$-functions revisited}
\author[]{J.KACZOROWSKI \lowercase{and} A.PERELLI}
\maketitle

{\bf Abstract.} The analytic properties of the standard twist $F(s,\alpha)$, where $F(s)$ belongs to a wide class of $L$-functions, are of prime importance in describing the structure of the Selberg class. In this paper we present a deeper study of such properties. In particular, we show that $F(s,\alpha)$ satisfies a functional equation of a new type, somewhat resembling that of the Hurwitz-Lerch zeta function. Moreover, we detect the finer polar structure of $F(s,\alpha)$, characterizing in two different ways the occurrence of finitely or infinitely many poles as well as giving a formula for their residues.



\vskip.5cm
\section{Introduction}

\smallskip
{\bf 1.1. The problem.}
Given an $L$-function $F(s)=\sum_{n\geq 1}a(n)n^{-s}$ of degree $d>0$ and a real number $\alpha>0$, let
\begin{equation}
\label{1-1}
F(s,\alpha)=\sum_{n=1}^{\infty}\frac{a(n)}{n^s} e(-\alpha n^{1/d}) \hskip2cm \big(e(x) = e^{2\pi ix}\big)
\end{equation}
 denote its standard twist. It is known that if $F(s)$ belongs to the extended Selberg class $\S^{\sharp}$, then $F(s,\alpha)$, initially defined for $\si>1$, extends to a meromorphic function on the whole complex plane $\mathbb C$. Moreover, there exists a discrete subset Spec$(F)$ of the positive real numbers, called the spectrum of $F(s)$, such that $F(s,\alpha)$ is entire if $\alpha\not\in$ Spec$(F)$, otherwise $F(s, \alpha)$ can have simple poles at the points (in case of normalized $F(s)$)
\begin{equation}
\label{1-2}
s_\ell = \frac{d+1}{2d} - \frac{\ell}{d} \hskip2cm \ell=0,1,2\dots
\end{equation}
 The first point $s_0$ is always a pole, the other may not. For these and other properties of $F(s,\alpha)$ the reader is referred to \cite{Ka-Pe/2005},\cite{Ka-Pe/2016a} and \cite{Ka-Pe/2016b}; see also Section 1.2 below  for exact definitions and a more explicit explanation of the notation.

\smallskip
The analytic properties of the standard twist are of prime importance in describing the structure of the Selberg class, see  \cite{Ka-Pe/1999a},\cite{Ka-Pe/2002b},\cite{Ka-Pe/2011}. A deeper study of such properties was initiated in \cite{Ka-Pe/half}, where the special case of the $L$-functions associated with half-integral weight modular forms was treated. Let $f$ be a cusp form of half-integral weight $\kappa=k/2$ with level $N$ and Fourier coefficients $a(n)$, where $k=2(h+1)+1$ with $h\in\{0,1,\dots\}$ and $4|N$, and let
\begin{equation}
\label{1-3}
F(s) = \sum_{n=1}^\infty \frac{a(n)}{n^{s+(\kappa-1)/2}}
\end{equation}
be its normalized Hecke $L$-function. Although $F(s)$ does not formally belong to $\S^\sharp$ in general, it enjoys very similar properties and in particular satisfies a Riemann-type functional equation. In \cite{Ka-Pe/half} we proved that the standard twist of $F(s)$, namely \eqref{1-1} with coefficients $a(n)n^{-(\kappa-1)/2}$ as in \eqref{1-3} and $d=2$, satisfies for $\alpha>0$ and $s\in\CC$ the functional equation
\begin{equation}
\label{1-4}
F(s,\alpha) = \frac{i^{-\kappa}}{\sqrt{2\pi}} \Big(\frac{\sqrt{N}}{4\pi}\Big)^{1-2s} \sum_{\ell=0}^h d_\ell \Gamma\big(2(1-s)-1/2- \ell\big) F^*_\ell(1-s,\alpha).
\end{equation}
Here $d_\ell$ are certain coefficients and $F^*_\ell(s,\alpha)$ certain entire functions closely related with $F(s)$. It follows in particular that for $\alpha\in$ Spec$(F)$, the standard twist $F(s,\alpha)$ has at most finitely many poles, located at the points $s_{\ell}$ with $\ell = 0,1,\dots, h$; see \cite{Ka-Pe/half}, \eqref{1-20}, \eqref{1-22} and Theorem 4 for details.
 
\smallskip
A closer analysis of the proof of \eqref{1-4} in \cite{Ka-Pe/half} reveals that these statements are consequences of the very special form of the functional equation of the $L$-function \eqref{1-3}. It  has a single $\Gamma$-factor of the form 
\[ 
\Gamma\Big(s+\frac{k-2}{4}\Big) \qquad \text{with} \ \ 2\nmid k,
\]
and the proof can be adapted to functional equations with a single $\Gamma$-factor of type
\begin{equation}
\label{1-5}
\Gamma\big(\frac{ds}{2}+\mu\big) \qquad \text{with} \ \ \mu= \frac{2(h+1)-d\pm1}{4} \geq 0,  \ h=0,1,\dots
\end{equation}
This is due to the fact that the argument is based on the explicit expression of the Mellin-Barnes integral
\begin{equation}
\label{1-6}
\frac{1}{2\pi i} \int_{(c)} \Gamma(\xi-w) \Gamma(w) \eta^{-w} \d w = \Gamma(\xi) (1+\eta)^{-\xi},
\end{equation}
where 
\begin{equation}
\label{1-7}
0<c<\Re (\xi) \quad \text{and} \quad |\arg \eta |<\pi;
\end{equation}
see (3.3.9) in Chapter 3 of Paris-Kaminski \cite{Pa-Ka/2001}. Nevertheless, one may speculate that results of type \eqref{1-4} should hold in general, but their proof requires some new ideas. In this paper we show that this hope is only partially justified, but fails in general. In particular, the situation changes drastically when we switch from half-integral weight cusp forms to the usual holomorphic cusp forms of integral weight $k$. The corresponding $L$-functions satisfy (after normalization) a functional equation with the single $\Gamma$-factor
 \[ 
 \Gamma\Big(s+\frac{k-1}{2}\Big),
 \]
but this seemingly unimportant difference has a major influence on the behavior of the corresponding standard twist. 
Very special cases of our Theorems 2 and 3 below imply that the standard twist $F(s,\alpha)$ corresponding to a holomorphic cusp form satisfies a certain functional equation, but  not as simple as \eqref{1-4}. Moreover, the polar structure changes dramatically. This time $F(s,\alpha)$ has always infinitely many poles for every $\alpha\in$ Spec$(F)$, as already observed in \cite{Kac/2004}.

\smallskip
The results in this paper are formulated in full generality, for arbitrary functions $F(s)$ of positive degree from $\S^{\sharp}$. In particular, we address and solve the following two problems:

(i) prove a functional equation for $F(s,\alpha)$, relating $s$ to $1-s$;

(ii) study the finer polar structure of $F(s,\alpha)$, in particular the existence of finitely or infinitely many poles at the points \eqref{1-2}.

\noindent
Our main result, see Theorem 2, is the solution of problem (i). Indeed, given $F\in\S^\sharp$ of degree $d$, we show that $F(s,\alpha)$ satisfies a functional equation of a new Hurwitz-Lerch type of degree $d$, of which \eqref{1-4} is a very special case with $d=2$. Moreover, in Section 1.3 we present a thorough discussion of the new functional equation, characterizing the occurrence of the special shape of type \eqref{1-4}, where the error function $H_k(s,\alpha)$ is not present, and deriving information on the polar structure of $F(s,\alpha)$. This leads to the solution of problem (ii).

\smallskip
Interestingly, despite the presence of an error function $H_k(s,\alpha)$ and other differences, such a functional equation still allows the deduction of many of the well known consequences of the functional equations of Riemann-type. In order to keep the size of the paper within reasonable limits, we postpone to a forthcoming paper a discussion of the applications to convexity bounds, distribution of zeros and other classical problems about $F(s,\alpha)$.

\smallskip
We conclude remarking that in the case of $F(s)$ of degree 1, thanks to the characterization of such functions in \cite{Ka-Pe/1999a}, $F(s,\alpha)$ reduces to a linear combination over certain Dirichlet polynomials of classical Hurwitz-Lerch zeta functions. Hence in this case the functional equation is immediately available, thus enabling a direct treatment of the above problems; see the recent paper by Zaghloul \cite{Zag/preprint}. Thus, although the general methods in this paper recover such results as special cases, it turns out once more that the degree 1 has rather special features inside the Selberg class.

\medskip
{\bf 1.2. Definitions and notation.}
Throughout the paper we write $s=\si+it$, an empty sum equals 0 and an empty product equals 1, logarithms have the principal value unless otherwise specified and $f(s)\equiv0$ means that $f(s)$ vanishes identically.

\smallskip
The extended Selberg class $\S^\sharp$ consists of non identically vanishing Dirichlet series $F(s)$, absolutely convergent for $\si>1$, such that $(s-1)^mF(s)$ is entire of finite order for some integer $m\geq0$ and satisfying a functional equation of type
\begin{equation}
\label{1-9}
F(s) \gamma(s) = \omega \overline{\gamma}(1-s) \overline{F}(1-s),
\end{equation}
where $\overline{f}(s)$ means $\overline{f(\overline{s})}$ and the $\gamma$-factor $\gamma(s)$ is defined by
\[
\gamma(s) = Q^s\prod_{j=1}^r\Gamma(\lambda_js+\mu_j) 
\]
with $|\omega|=1$, $Q>0$, $r\geq0$, $\lambda_j>0$ and $\Re{\mu_j}\geq0$. We refer to Selberg \cite{Sel/1989}, Conrey-Ghosh \cite{Co-Gh/1993} and to our survey papers \cite{Ka-Pe/1999b},\cite{Kac/2006},\cite{Per/2005},\cite{Per/2004},\cite{Per/2010},\cite{Per/2017} for definitions, examples and the basic theory of the Selberg class. We recall that degree $d$, conductor $q$, $\xi$-invariant $\xi_F$ and root number $\omega_F$ of $F(s)$ are the invariants defined by
\begin{equation}
\label{1-10}
\begin{split}
d=d_F &:=2\sum_{j=1}^r\lambda_j, \quad q=q_F:= (2\pi)^dQ^2\prod_{j=1}^r\lambda_j^{2\lambda_j}, \\
 \quad \xi_F = 2\sum_{j=1}^r &(\mu_j-\frac12) = \eta_F+id \theta_F,  \quad \omega_F=\omega\prod_{j=1}^r \lambda_j^{-2i\Im(\mu_j)}
\end{split}
\end{equation}
with $\eta_F,\theta_F\in\RR$. The invariant $\theta_F$ is called the internal shift of $F(s)$, and the classical $L$-functions have $\theta_F=0$. Moreover, $m_F$ is the order of pole of $F(s)$ at $s=1$ and 
\begin{equation}
\label{1-11}
 \sum_{m=1}^{m_F} \frac{\gamma_m}{(s-1)^m} 
\end{equation}
is its polar part. In this paper {\it we always assume that the degree $d$ is positive}, and hence $d\geq1$ thanks to the results in \cite{Co-Gh/1993} and \cite{Ka-Pe/1999a}.

\smallskip
The spectrum of $F(s)$ is defined as
\begin{equation}
\label{1-12}
\text{Spec}(F) = \{\alpha>0: a(n_\alpha)\neq0\} = \Big\{d\big(\frac{m}{q}\big)^{1/d}: m\in\NN \ \text{with} \ a(m)\neq 0\Big\},
\end{equation}
where $a(n)$ are the coefficients of $F(s)$ and
\begin{equation}
\label{1-13}
n_\alpha = q d^{-d} \alpha^{d}, \quad a(n_\alpha)=0 \ \text{if} \ n_\alpha\not\in \NN.
\end{equation}
We recall that the standard twist $F(s,\alpha)$ is entire if $\alpha\not\in$ Spec$(F)$, while for $\alpha\in$ Spec$(F)$ it is meromorphic on $\CC$ with at most simple poles at the points 
\begin{equation}
\label{1-14}
s^*_\ell = s_\ell - i \theta_F,
\end{equation}
where $s_\ell$ and $\theta_F$ are as in \eqref{1-2} and \eqref{1-10}, respectively. The residue of $F(s,\alpha)$ at $s=s^*_\ell$ is denoted by $\rho_\ell(\alpha)$. It is know that $\rho_0(\alpha)\neq0$ when $\alpha \in$ Spec$(F)$. Moreover, $F(s,\alpha)$ has polynomial growth on every vertical strip, although the known bounds are weak in general. We refer to \cite{Ka-Pe/2005} and \cite{Ka-Pe/2016a} for these and other results on $F(s,\alpha)$.

\smallskip
Let
\begin{equation}
\label{1-15}
S_F(s) := 2^r \prod_{j=1}^r \sin(\pi(\lambda_js + \mu_j)) = \sum_{j=-N}^N a_j e^{i\pi d_F\omega_j s}
\end{equation}
with a certain $N\in\NN$, $a_j\in\CC$ and
\begin{equation}
\label{1-16}
-\frac12 =\omega_{-N} < \dots <\omega_N = \frac12
\end{equation}
satisfying $\omega_{-j}=-\omega_j$. Moreover $a_j\neq0$ for $j\neq0$, and $N\geq1$ since $r\geq1$. We also define
\begin{equation}
\label{1-17}
h_F(s) = \frac{\omega}{(2\pi)^r} Q^{1-2s} \prod_{j=1}^r \big(\Gamma(\lambda_j(1-s)+\overline{\mu}_j) \Gamma(1-\lambda_js-\mu_j)\big).
\end{equation}
It turns out, see \eqref{1-22} and Theorem 2 of \cite{Ka-Pe/2000}, that {\it the functions} $S_F(s)$ {\it and} $h_F(s)$ {\it are invariants}. Moreover, thanks to the reflection formula of the $\Gamma$ function, the functional equation \eqref{1-9} of $F(s)$ can be written in the following {\it invariant form} 
\begin{equation}
\label{1-18}
F(s) = h_F(s) S_F(s) \overline{F}(1-s).
\end{equation}

\smallskip
{\bf Remark 1.1.} Write \eqref{1-18} as
\[
 S_F(s) \overline{F}(1-s) = \frac{F(s)}{h_F(s)}
\]
and observe that the left hand side is certainly holomorphic for $s\neq0$, and so is the right hand side for $s\neq 1$. Hence, in particular, it follows that $S_F(s) \overline{F}(1-s)$ {\it is entire}. \qed

\smallskip
Finally, with $\gamma_m$ as in \eqref{1-11}, we consider the {\it residual function}
\begin{equation}
\label{1-19}
R(s,\alpha) =  (2\pi i \alpha)^{-ds} \sum_{m=1}^{m_F} d^m\gamma_m \sum_{h=0}^{m-1} \frac{(-1)^h \log^h (2\pi i\alpha) \Gamma^{(m-1-h)}(ds)}{h! (m-1-h)!},
\end{equation}
hence $R(s,\alpha)\equiv0$ if $F(s)$ is entire, and for $\ell=0,1,\dots$ let
\begin{equation}
\label{1-20}
\overline{F}_\ell(s,\alpha) = \sum_{j=-N}^N a_j e^{i\pi d\omega_j (1-s)} \sideset{}{^\flat} \sum_{n\geq1} \frac{\overline{a(n)}}{n^s} \left(1+ e^{i\pi(\frac12-\omega_j)} \left(\frac{n_\alpha}{n}\right)^{1/d} \right)^{d(1-s-s_\ell^*)},
\end{equation}
where the symbol $^\flat$ in the inner sum indicates that if $j=-N$ then the term $n=n_\alpha$ is omitted. Hence $\overline{F}_\ell(s,\alpha)$ is well defined since $1+ e^{i\pi(\frac12-\omega_j)} \left(\frac{n_\alpha}{n}\right)^{1/d} \neq 0$ always. Note that if $\alpha \not\in$ Spec$(F)$ we may omit $^\flat$, since $a(n_\alpha)=0$ in this case. Note also that the inner sum in \eqref{1-20} is a general Dirichlet series with complex frequencies, absolutely convergent for $\si>1$, and
\begin{equation}
\label{1-21}
 \left(1+ e^{i\pi(\frac12-\omega_j)} \left(\frac{n_\alpha}{n}\right)^{1/d} \right)^{d(1-s-s_\ell^*)} = e^{d(1-s-s_\ell^*) \Log\big(1+ e^{i\pi(\frac12-\omega_j)} \left(\frac{n_\alpha}{n}\right)^{1/d} \big)}
\end{equation}
where the branch of Log on $\overline{\HH} \setminus\{0\}$ has argument in $[0,\pi]$, $\HH$ being the upper half-plane. 

\medskip
{\bf 1.3. Main results.}
We start with the properties of the functions $\overline{F}_\ell(s,\alpha)$ defined in \eqref{1-20}. As one can guess from the definition, after expanding the left hand side of \eqref{1-21} these functions are close to suitable ``stratifications'' of $\overline{F}(s)$; see \eqref{2-1} and \eqref{2-4}.

\medskip
{\bf Theorem 1.} {\sl Let $F\in\S^\sharp$ with degree $d\geq1$, $\alpha>0$ and $\ell=0,1,2,\dots$ Then $\overline{F}_\ell(s,\alpha)$ is an entire function, not identically vanishing. Moreover, uniformly for $\si$ in any bounded interval, as $|t|\to\infty$ we have
\[
\overline{F}_\ell(s,\alpha) \ll e^{\frac{\pi}{2} d |t|} |t|^{c(\si)} 
\]
with a certain $c(\si)\geq0$ independent of $\ell$ and $\alpha$, satisfying $c(\si)=0$ for $\si>1$.}

\medskip
The dependence of the implicit constant on $\alpha$ can be made explicit, starting with an explicit bound in $\delta$ in Lemma A, see the Appendix, and then following the proof of Theorem 1 in Section 2.2. The dependence on $\ell$ and bounds for $c(\si)$ can also be made explicit, again following the arguments in Section 2.2.

\smallskip
Before stating the functional equation of the standard twist $F(s,\alpha)$, we introduce a new set of invariants denoted by $d_\ell$, related to an asymptotic expansion of the function $h_F(s)$ in \eqref{1-17}. Indeed, given $F\in\S^\sharp$ with degree $d$ and conductor $q$, for large $|s|$ outside an arbitrarily small angular region containing the positive real axis we have
\begin{equation}
\label{1-22}
h_F(s) \approx \frac{\omega_F}{\sqrt{2\pi}} \left(\frac{q^{1/d}}{2\pi d}\right)^{d(\frac12-s)} \sum_{\ell=0}^\infty d_\ell \Gamma\big(d(s_\ell^*-s)\big),
\end{equation}
where $s_\ell$ are defined in \eqref{1-2} and $\approx$ means that cutting the sum at $\ell=M$ one gets a meromorphic remainder which is $\ll$ than the modulus of the $M$-th term times $1/|s|$. We shall prove a suitable version of \eqref{1-22} in Section 3.2 when $\theta_F=0$, see \eqref{3-16} and Remark 3.1; see also Section 3.6 for the general case. Although obtained by an application of the Stirling formula, such an expansion is somehow non-standard, and is crucial for our method. Indeed, in Section 3 it will allow us to treat the general case of $F\in\S^\sharp$ by the Mellin-Barnes integral in \eqref{1-6}, which is the heart of the treatment of the special case in \cite{Ka-Pe/half}. In view of their relevance for the structure of the functional equation of $F(s,\alpha)$, the invariants $d_\ell=d_\ell(F)$ are called the {\it structural coefficients} of $F(s)$. In particular we have, see \eqref{3-64}, that
\begin{equation}
\label{1-23}
d_0= d^{id\theta_F}.
\end{equation}

\smallskip
We are now ready to state the functional equation of $F(s,\alpha)$ for any $F\in\S^\sharp$, thus providing a vast extension of \eqref{1-4}. We use the notation introduced in Sections 1.1 and 1.2, we write
\begin{equation}
\label{1-24}
e_\ell = \frac{d-1}{2} +\ell +i\theta_F
\end{equation}
and we recall that the residual function $R(s,\alpha)$ in \eqref{1-19} is identically vanishing if $F(s)$ is entire.

\medskip
{\bf Theorem 2.} {\sl Let $F\in\S^\sharp$ with $d\geq1$ and let $\alpha>0$. Then for any integer $k\geq 0$ and $s$ in the strip $s_{k+1} < \si < s_k$ we have
\begin{equation}
\label{1-26}
\begin{split}
F(s,\alpha) =   \frac{\omega_F}{\sqrt{2\pi}} \left(\frac{q^{1/d}}{2\pi d}\right)^{d(\frac12-s)} \sum_{\ell=0}^k &d_\ell \Gamma\big(d(1-s) - e_\ell \big) \overline{F}_{\ell}(1-s,\alpha) \\
&+ R(1-s,\alpha) + H_k(s,\alpha),
\end{split}
\end{equation}
where the coefficients $d_\ell$ are defined by \eqref{1-22} and the function $H_k(s,\alpha)$ is  holomorphic in the above strip and meromorphic over $\CC$ with all poles in a horizontal strip of bounded height. Moreover, there exists $\theta =\theta(d)>0$ such that for any $\si\in[s_{k+1},s_k]\cap (-\infty,0)$ we have
\[
H_k(s,\alpha) \ll |t|^{-\theta} \hskip1.5cm \text{as} \ \ |t|\to\infty.
\]}

\medskip
In view of the error function $H_k(s,\alpha)$, we may regard \eqref{1-26} as an approximate form of a general Hurwitz-Lerch functional equation of degree $d$. The constant in the $\ll$-symbol may depend on $\si$, $F(s)$ and $\alpha$. The value of $\theta$ is usually close to $1/2$ but may become smaller when $\si$ is near 0; we refer to Lemma 3.2 and Remark 3.4 in Section 3.5 for explicit bounds for $H_k(s,\alpha)$. Finally, by \eqref{1-2} and \eqref{1-24} the $\Gamma$-factors in \eqref{1-26} may also be written as $\Gamma(d(s_\ell^*-s))$.

\smallskip
From the proof of Theorem 2 we obtain the following explicit expression for the residues $\rho_\ell(\alpha)$ of $F(s,\alpha)$ at the potential poles $s_\ell$ in \eqref{1-2}, when $\alpha\in$ Spec$(F)$.

\medskip
{\bf Theorem 3.} {\sl Let $F\in\S^\sharp$ with $d\geq 1$, $\alpha\in$ {\rm Spec}$(F)$, $\ell=0,1,\dots$ and $d_\ell$ be defined by \eqref{1-22}. Then
\[
\rho_\ell(\alpha) = \frac{d_\ell}{d} \frac{\omega_F}{\sqrt{2\pi}} e^{-i\frac{\pi}{2}(\xi_F+ds_\ell^*)} \left(\frac{q^{1/d}}{2\pi d}\right)^{\frac{d}{2}-ds_\ell^*} \frac{\overline{a(n_\alpha)}}{n_\alpha^{1-s_\ell^*}}.
\]
In particular,  the set of poles of $F(s,\alpha)$ is independent of $\alpha$ and equals $\{s_\ell^*: d_\ell\neq 0\}$.}

\medskip
We have seen in Section 1.1 that there are cases where the function $H_k(s,\alpha)$ in \eqref{1-26} vanishes identically for some value of $k$; see \eqref{1-4}. Now we study the conditions under which this phenomenon occurs; for simplicity we consider only the case where $\theta_F=0$. Suppose that there exist $\alpha_0>0$ and a {\it minimal} integer $h\geq 0$ such that
\begin{equation}
\label{1-27}
H_k(s,\alpha_0) \equiv 0 \ \text{{\it for every $k\geq h$}}.
\end{equation}
Then, since the functions $\overline{F}_\ell(1-s,\alpha_0)$ are not identically vanishing by Theorem 1, it follows from \eqref{1-26} that
\begin{equation}
\label{1-28}
\text{{\it $d_\ell=0$ for every $\ell\geq h+1$.}}
\end{equation}
Moreover, $H_k(s,\alpha_0)$ can be explicitly computed for $0\leq k \leq h-1$,  see \eqref{4-18}. Since the coefficients $d_\ell$ are independent of $\alpha$, assuming \eqref{1-27} we get from Theorem 3 and \eqref{1-28} that
\begin{equation}
\label{1-29}
\text{{\it  the poles of $F(s,\alpha)$ are at the points $\{s_\ell^*: 0\leq \ell\leq h, d_\ell\neq0\}$ for every $\alpha\in$ {\rm Spec}$(F)$}}.
\end{equation}
It turns out that the opposite implication holds as well, even in the strong form that assuming \eqref{1-29} then $H_k(s,\alpha)\equiv0$ for all $k\geq h$ and every $\alpha>0$. Hence we give the following definition.

\medskip
{\bf Definition 1.1.} Let $F\in\S^\sharp$ with $d\geq 1$ and $\theta_F=0$. We say that $F(s,\alpha)$ satisfies a {\it strict functional equation} if there exists an integer $h\geq 0$, minimal and independent of $\alpha$, such that $H_k(s,\alpha)\equiv0$ for every $k\geq h$ and all $\alpha>0$. \qed

\medskip
\noindent
Thus, in view of \eqref{1-26} the strict functional equation of $F(s,\alpha)$ has the form
\begin{equation}
\label{1-30}
F(s,\alpha) =   \frac{\omega_F}{\sqrt{2\pi}} \left(\frac{q^{1/d}}{2\pi d}\right)^{d(\frac12-s)} \sum_{\ell=0}^{h} d_\ell \Gamma\big(d(1-s) - e_\ell \big) \overline{F}_{\ell}(1-s,\alpha) + R(1-s,\alpha)
\end{equation}
for every $s\in\CC$ and $\alpha>0$. We refer to Section 4.3 for several remarks about it. 

\smallskip
The next result gives two characterizations of the occurrence of the strict functional equation, one in terms of the poles of $F(s,\alpha)$ and the other in terms of the $\gamma$-factors of $F(s)$, but first we need to introduce a further definition. Let $N\geq1$ and $n_j\geq0$, $j=1,\dots,N$, be integers; we say that $n_1,\dots,n_N$ form a {\it compatible system} if
\begin{equation}
\label{1-31}
\text{ (a) \ $n_i\not\equiv n_j$ {\rm (mod $2N$)} for every $i\neq j$,} \quad \text{ (b) \ $n_i\not\equiv 1-n_j $ {\rm (mod $2N$)} for every $i,j$.}
\end{equation}
Conditions (a) and (b) imply that the sets
\[
\{n_j \ \text{(mod $2N$)}\} \quad \text{and} \quad \{2N+1-n_j \ \text{(mod $2N$)}\}
\]
are disjoint and their union is $\{0,1,\dots,2N-1 \ \text{(mod $2N$)} \}$. Thus there exists $n_0\in\ZZ$ such that
\begin{equation}
\label{1-32}
\{n_j, 2N+1-n_j \ \text{with} \ j=1,\dots,N\} = \{n_0+j + 2N\nu_j, j=0,\dots,2N-1\},
\end{equation}
where the integers $\nu_j\geq0$ are such that
\begin{equation}
\label{1-33}
N_0= \sum_{j=0}^{2N-1} \nu_j
\end{equation}
is minimal. The integer $N_0\geq0$ has an interesting property, as shown in statement (iii) below.

\medskip
{\bf Theorem 4.} {\sl Let $F\in\S^\sharp$ with $d\geq1$ and $\theta_F=0$, and let $\alpha>0$. Then the following statements are equivalent.

(i) $F(s,\alpha)$ satisfies a strict functional equation of the form \eqref{1-30};

(ii) for every $\alpha\in$ {\rm Spec}$(F)$ all the poles of $F(s,\alpha)$ are at the points $\{s_\ell: 0\leq \ell\leq h, d_\ell\neq0\}$ and $d_h\neq0$;

(iii) $F(s)$ has a $\gamma$-factor of the form
\[
\gamma(s) = Q^s \prod_{j=1}^N \Gamma\Big(\frac{d}{2N}s + \frac{2n_j - d-1}{4N}\Big),
\]
where $Q>0$, $N\geq1$ and the integers $n_j$ satisfy $n_j\geq (d+1)/2$ and form a compatible system with $N_0=h$.}

\medskip
Thus, in particular, Theorem 4 gives a characterization of the occurrence of the strict functional equation and of the polar structure of $F(s,\alpha)$ in terms of the $\gamma$-factors of $F(s)$. Moreover, in view of (ii) and Theorem 3, a third statement equivalent to (i) is that the structural coefficients $d_\ell$ satisfy $d_h\neq 0$ and $d_\ell =0$ for every $\ell\geq h+1$. We conclude with two remarks.

\smallskip
{\bf Remark 1.2.} In view of (iii), \eqref{1-32} and \eqref{1-33} provide an algorithm to compute the value of $h$ in \eqref{1-30} and in (ii), starting from a suitable $\gamma$-factor of $F(s)$. We refer to \cite{Ka-Pe/2000} for the transformation techniques of $\gamma$-factors, which can help deciding if a function $F(s)$ has a $\gamma$-factor of type (iii). Note that \eqref{1-5} and the $\gamma$-factor of the function $F(s)$ in example \eqref{1-3} are of the form (iii), and we refer to Remark 7 in Section 4.3 for an implementation of the algorithm in these cases. \qed

\smallskip
{\bf Remark 1.3.}  At present we don't know examples of $L$-functions whose standard twist satisfies a strict functional equation, other than the degree 1 $L$-functions and the Hecke $L$-functions of half-integral weight cusp forms in \eqref{1-3}. In particular, we don't know if there exist $L$-functions having a $\gamma$-factor of type \eqref{1-5} with integer degree $d\geq3$. It would be of course interesting to exhibit an example of this type, or to prove non-existence; the second task appears to be difficult at present. \qed

\medskip
{\bf Acknowledgements.} 
This research was partially supported by the Istituto Nazionale di Alta Matematica, by the MIUR grant PRIN-2015 {\sl ``Number Theory and Arithmetic Geometry''} and by grant 2017/25/B/ST1/00208 {\sl ``Analytic methods in number theory''}  from the National Science Centre, Poland.

\bigskip
\section{Proof of Theorem 1}

\smallskip
We prove Theorem 1 in two steps.

\medskip
{\bf 2.1. Analytic continuation and nonvanishing.}
We split the sum over $n$ in the definition \eqref{1-20} of $\overline{F}_\ell(s,\alpha)$ into $n\leq n_\alpha$ and $n>n_\alpha$, thus getting that
\begin{equation}
\label{2-1}
\overline{F}_\ell(s,\alpha) = \overline{F}_\ell^{(1)}(s,\alpha) + \overline{F}_\ell^{(2)}(s,\alpha),
\end{equation}
say. Clearly, $\overline{F}_\ell^{(1)}(s,\alpha)$ is an entire function of $s$. 

\smallskip
The treatment of $\overline{F}_\ell^{(2)}(s,\alpha)$ is based on Lemma A in the Appendix, with the choices
\begin{equation}
\label{2-2}
\rho = d(1-s-s_\ell^*) \quad \text{and} \quad z= z_{n,j}:=e^{i\pi(\frac12-\omega_j)} \left(\frac{n_\alpha}{n}\right)^{1/d};
\end{equation}
note that $|z| \leq 1-\delta$ for some $\delta=\delta(\alpha)>0$ since $n>n_\alpha$. Assuming $\si>1$, Lemma A implies that for any integer $R\geq0$ we have
\begin{equation}
\label{2-3}
\begin{split}
\overline{F}_\ell^{(2)}(s,\alpha) &=  \sum_{j=-N}^N a_j e^{i\pi d\omega_j (1-s)} \sum_{r=0}^R {\rho \choose r} e^{i\pi r(\frac12-\omega_j)} n_\alpha^{r/d}  \sum_{n>n_\alpha} \frac{\overline{a(n})}{n^{s+r/d}}  \\
&\hskip2cm + \sum_{j=-N}^N a_j e^{i\pi d\omega_j (1-s)} \sum_{n>n_\alpha}  \frac{\overline{a(n})}{n^s} Q_R(z_{n,j},\rho)\\
&=  \sum_{j=-N}^N a_j e^{i\pi d\omega_j (1-s)} \sum_{r=0}^R {\rho \choose r} e^{i\pi r(\frac12-\omega_j)} n_\alpha^{r/d} \Big(\overline{F}\big(s+\frac{r}{d}\big) - \sum_{n\leq n_\alpha} \frac{\overline{a(n})}{n^{s+r/d}} \Big) \\
&\hskip 2cm + \sum_{j=-N}^N a_j e^{i\pi d\omega_j (1-s)} \sum_{n>n_\alpha}  \frac{\overline{a(n})}{n^s} Q_R(z_{n,j},\rho).
\end{split}
\end{equation}
Moreover, since $\max_{0\leq u \leq 1} \arg(1+uz_{n,j}) \to 0$ as $n\to\infty$, by \eqref{6-1} in the Appendix the last sum is absolutely convergent for $\si> 1-(R+1)/d$ and hence is holomorphic on this half-plane. Therefore, by \eqref{1-15}, for $\si>1-(R+1)/d$ we get
\begin{equation}
\label{2-4}
\begin{split}
\overline{F}_\ell^{(2)}(s,\alpha) &=  \sum_{r=0}^R {\rho \choose r} e^{i\pi r/2} n_\alpha^{r/d}  \sum_{j=-N}^N a_j e^{i\pi d\omega_j (1-(s+\frac{r}{d}))} \overline{F}\big(s+\frac{r}{d}\big) + h(s) \\
&=  \sum_{r=0}^R {\rho \choose r} e^{i\pi r/2} n_\alpha^{r/d} S_F\big(1-(s+\frac{r}{d})\big) \overline{F}\big(s+\frac{r}{d}\big) + h(s)
\end{split}
\end{equation}
with a certain holomorphic function $h(s)$. But $S_F(1-s) \overline{F}(s)$ is entire, see Remark 1.1; since $R$ is arbitrary, from \eqref{2-1} and \eqref{2-4} we see that $\overline{F}_\ell(s,\alpha)$ is entire as well.

\smallskip
Now we show that $\overline{F}_\ell(s,\alpha)$ is not identically vanishing. Recalling \eqref{1-21} and \eqref{2-2} and writing
\[
r_{n,j} = n^{1/d} |1+z_{n,j}|, \quad \theta_{n,j}= \arg(1+ z_{n,j}),
\]
the generic term of $\overline{F}_\ell(\si,\alpha)$ has the form
\begin{equation}
\label{2-5}
a_j \overline{a(n)} A_{n,j} \big(r_{n,j}e^{i(\pi\omega_j+\theta_{n,j})}\big)^{-d\si},
\end{equation}
where $A_{n,j}\in\CC\setminus\{0\}$ does not depend on $\si$. Thanks to the triangle inequality, for $j\neq-N$ and $1\leq n\neq n_\alpha$ we have
\[
r_{n,j} > n^{1/d}  \Big|1-\Big(\frac{n_\alpha}{n}\Big)^{1/d}\Big| = r_{n,-N}.
\]
Hence we define $1\leq n_0 \neq n_\alpha$ to be an integer such that 
\begin{equation}
\label{2-6}
n_0^{1/d} \Big|1-\Big(\frac{n_\alpha}{n_0}\Big)^{1/d}\Big| = \min_{\substack{n\neq n_\alpha \\ a(n)\neq 0}} \Big\{n^{1/d} \Big|1-\Big(\frac{n_\alpha}{n}\Big)^{1/d}\Big| \Big\}.
\end{equation}
Clearly, the expression in \eqref{2-6} is $>0$ and $a(n_0)\neq 0$. Such $n_0$ exists since the expression in curly brackets tends to $\infty$ as $n\to\infty$. Moreover, there are at most two integers $n_0$ satisfying \eqref{2-6}, since the function $\xi|1-A/\xi|$ with $A>0$ has only one minimum for $\xi>0$, namely $\xi=A$. Clearly, if there are two such integers $n_1<n_2$, then $n_\alpha\in(n_1,n_2)$. Therefore, for every $n\neq n_\alpha$ and $-N\leq j\leq N$ we have that
\begin{equation}
\label{2-7}
\text{$r_{n,j} > r_0$  if $(n,j)\neq (n_0,-N)$ and, by definition, $r_{n_0,-N}=r_0$},
\end{equation}
where $n_0$ denotes the generic solution of \eqref{2-6}. In view of the symbol $^\flat$ in \eqref{1-20}, if $n=n_\alpha$ and $a(n_\alpha)\neq0$ we consider only $j\neq -N$. In this case by   \eqref{1-16} we have $0\leq \pi(\frac{1}{2} -\omega_j) < \pi$, thus thanks to triangle inequality we see that
\begin{equation}
\label{2-8}
\rho_0 :=\min_{\substack{j\neq -N \\ a_j\neq 0}} r_{n_\alpha,j} = r_{n_\alpha,j_0}, \ \ \text{where} \ \ j_0=
\begin{cases}
-N+1 & \text{if} \ N\geq 2 \\
1 & \text{if} \ N=1 \ \text{and} \ a_0=0,
\end{cases}
\end{equation}
and the minimum is attained only at $j=j_0$.

\smallskip
From \eqref{2-7} and \eqref{2-8} we deduce that the dominating terms in $\overline{F}_\ell(\si,\alpha)$ correspond to the terms \eqref{2-5} with $(n,j) = (n_1,-N), (n_2,-N), (n_\alpha,j_0)$, meaning that as $\si\to+\infty$
\begin{equation}
\label{2-9}
\begin{split}
\overline{F}_\ell(\si,\alpha) &= A(r_0e^{i\frac{\pi}{2}})^{-d\si} + B(r_0e^{-i\frac{\pi}{2}})^{-d\si} + C(\rho_0 e^{i\lambda})^{-d\si}  + o(r_0^{-d\si}) + o(\rho_0^{-d\si}) \\
&= M(\si) +o(r_0^{-d\si}) + o(\rho_0^{-d\si}),
\end{split}
\end{equation}
say, with certain complex coefficients $A,B,C$ independent of $\si$ and not all zero (actually, at least one among $A$ and $B$ is non-zero), and a certain $\lambda\in\RR$. Indeed, the first main term in \eqref{2-9}, if it exists, corresponds to $n_1<n_\alpha$ and $j=-N$, hence $\omega_{-N}=-1/2$ and $\theta_{n_1,-N}=\pi$. The second one, if it exists, corresponds to  $n_2>n_\alpha$ and $j=-N$, hence $\omega_{-N}=-1/2$ and $\theta_{n_2,-N}=0$, while the third one corresponds to $n=n_\alpha$ and $j=j_0$, and is present if and only if $\alpha\in$ Spec$(F)$. Clearly, at least one of the first two terms exists and has non-zero coefficient. Moreover, the sum of all the other terms in the absolutely convergent series $\overline{F}_\ell(\si,\alpha)$ is obviously $o(r_0^{-d\si})$ or $o(\rho_0^{-d\si})$; hence \eqref{2-9} follows. 

\smallskip
In view of \eqref{2-9}, if $\alpha\in$ Spec$(F)$ and $\rho_0< r_0$ then clearly $|M(\si)| \sim |C| \rho_0^{-d\si}$, while if $\rho_0> r_0$ then $|M(\si)| =\Omega(r_0^{-d\si})$. This is indeed obvious if $A$ or $B$ vanishes, while if $AB\neq0$ then 
\begin{equation}
\label{2-10}
P(i\si) := Ae^{-\frac{\pi}{2}di\si} + Be^{\frac{\pi}{2}di\si}
\end{equation}
is the restriction to the imaginary axis of the general Dirichlet polynomial $P(u+i\si)$ with real frequencies. Thus $P(u+i\si)$ is almost periodic and our assertion follows. The case $\alpha\not\in$ Spec$(F)$ is similar. Therefore, $\overline{F}_\ell(s,\alpha)$ is not identically vanishing in these cases. If $\alpha\in$ Spec$(F)$ and $\rho_0=r_0$ we need a closer inspection of the third term in $M(\si)$. We first note that if $j_0=1$ in \eqref{2-8} then $\rho_0=2n_\alpha^{1/d}$ and $\lambda=\pi/2$, while $A=0$ since otherwise $r_0=n_\alpha^{1/d}-n_1^{1/d}<\rho_0$; hence $M(\si)$ is similar to \eqref{2-10} and our assertion follows in this case. Therefore, if $\rho_0=r_0$ then $-1/2 < \omega_{j_0} < 1/2$ and hence $0<\pi(\half -\omega_{j_0}) < \pi$. Writing
\[
\pi(\half -\omega_{j_0}) = \pi(1-\delta) \quad \text{with} \ 0<\delta<1,
\]
by elementary geometry we have that 
\[
\theta_{n_\alpha,j_0} = \arg(1+e^{i\pi(1-\delta)}) = \frac{\pi}{2}(1-\delta)
\]
and hence by \eqref{2-5}
\[
\lambda = \pi\omega_{j_0} + \theta_{n_\alpha,j_0} = \pi(\delta-\half) +  \frac{\pi}{2}(1-\delta) = \frac{\pi}{2}\delta.
\]
Therefore, when $\rho_0=r_0$ the three (or two) exponentials in $M(\si)$ are all distinct, thus, arguing as before, $|M(\si)|=\Omega(r_0^{-d\si})$ in this case as well, and our assertion on the nonvanishing of $\overline{F}_\ell(s,\alpha)$ follows. The first part of Theorem 1 is therefore proved.

\medskip
{\bf 2.2. Estimates.}
We first estimate the term in the last row of \eqref{2-3}. It is obvious that if $z\in\overline{\HH}\setminus\{0\}$ and $0\leq u\leq 1$ then
\[
0 \leq \arg(1+uz) \leq \arg(z) \leq \pi
\]
and hence, recalling \eqref{2-2}, for $0\leq u\leq 1$ we have that
\begin{equation}
\label{2-11}
0\leq \arg(1+uz_{n,j}) \leq \pi(\half-\omega_j).
\end{equation}
Therefore, by Lemma A in the Appendix with the choices in \eqref{2-2}, for $t\geq0$ we obtain 
\[
Q_R(z_{n,j},\rho) \ll \frac{(1+t)^{R+1}}{n^{\frac{R+1}{d}}} e^{dt\pi(\frac{1}{2}-\omega_j)},
\]
uniformly for $\si$ in any bounded interval. Thus for $t\geq0$ and any non-negative  integer $R$ satisfying $R>d-1-d\si$ we get
\[
\sum_{j=-N}^N a_j e^{i\pi d\omega_j (1-s)} \sum_{n>n_\alpha}  \frac{\overline{a(n})}{n^s} Q_R(z_{n,j},\rho) \ll e^{\frac{\pi}{2}dt} (1+t)^{R+1} \sum_{n>n_\alpha} \frac{|a(n)|}{n^{\si+\frac{R+1}{d}}} \ll e^{\frac{\pi}{2}dt} (1+t)^{R+1}
\]
uniformly for $\si$ in any bounded interval. Similarly, for $t\leq0$ from Lemma A and \eqref{2-11} we have
\[
Q_R(z_{n,j},\rho) \ll \frac{(1+|t|)^{R+1}}{n^{\frac{R+1}{d}}},
\]
hence in view of \eqref{1-16} for $R>d-1-d\si$ we get again
\[
\sum_{j=-N}^N a_j e^{i\pi d\omega_j (1-s)} \sum_{n>n_\alpha}  \frac{\overline{a(n})}{n^s} Q_R(z_{n,j},\rho) \ll e^{\frac{\pi}{2}d|t|} (1+|t|)^{R+1}.
\]
The contribution to $\overline{F}_\ell^{(2)}(s,\alpha)$ of the term in the last row of \eqref{2-3} is therefore
\begin{equation}
\label{2-12}
\ll e^{\frac{\pi}{2}d|t|} (1+|t|)^{R+1}
\end{equation}
uniformly for $\si$ in any bounded interval, where $R\geq0$ is any  integer satisfying $R>d-1-d\si$.

\smallskip
The treatment of $\overline{F}_\ell^{(1)}(s,\alpha)$ in \eqref{2-1} is similar. Indeed, thanks to \eqref{2-11} with $u=1$ we may write $\arg(1+z_{n,j}) = \pi(\half-\omega_j) - \eta_{n,j}$ with some $0\leq \eta_{n,j} \leq \pi$. Hence, without using Lemma A, we get
\begin{equation}
\label{2-13}
\overline{F}_\ell^{(1)}(s,\alpha) \ll \sum_{j=-N}^N e^{\pi d\omega_jt} \sideset{}{^\flat} \sum_{n\leq n_\alpha} e^{(\pi(\frac{1}{2} -\omega_j)-\eta_{n,j})dt} \ll  \sum_{j=-N}^N \sum_{n\leq n_\alpha} e^{(\frac{\pi}{2}-\eta_{n,j})dt} \ll e^{\frac{\pi}{2}d|t|},
\end{equation}
uniformly for $\si$ in any bounded interval.

\smallskip
Finally, the contribution to $\overline{F}_\ell^{(2)}(s,\alpha)$ of the term in the third row of \eqref{2-3} is
\begin{equation}
\label{2-14}
\begin{split}
\ll \sum_{j=-N}^N e^{\pi d\omega_j t} \sum_{r=0}^R (1+|t|)^r \Big( \big|F\big(s+\frac{r}{d}\big)\big| +  1 \Big) \ll e^{\frac{\pi}{2}|t|} \max_{0\leq r \leq R} (1+|t|)^{r+\mu(\si +\frac{r}{d}) +\epsilon}
\end{split}
\end{equation}
uniformly for $\si$ in any bounded interval, where $\mu(\si)$ is the Lindel\"of $\mu$-function of $F(s)$ and $\epsilon>0$ is arbitrarily small.

\smallskip
Now we choose 
\[
R=\max(0,[d(1-\si)]),
\]
so that $R>d-1-d\si$ and hence \eqref{2-12} holds. From \eqref{2-1},\eqref{2-3} and \eqref{2-12}-\eqref{2-14} we therefore have that the bound in Theorem 1 is satisfied with the choice
\[
c(\si) = \max \Big(R+1, \max_{0\leq r\leq R} \big\{r+\mu\big(\si+\frac{r}{d}\big)+\epsilon\big\}\Big).
\]
Moreover, when $\si>1$ the same argument leading to \eqref{2-13} shows that
\[
\overline{F}_\ell(s,\alpha) \ll e^{\frac{\pi}{2}d|t|},
\]
hence $c(\si)=0$ for $\si>1$ and Theorem 1 follows.

\bigskip
\section{Proof of Theorems 2 and 3}

\smallskip
For simplicity we first assume that $\theta_F=0$, which is usually satisfied by the classical $L$-functions. Then, in Section 3.6, we show how the general case follows from this special case.

\medskip
{\bf 3.1. Set up of Theorem 2.} 
As customary in the study of the standard twist, we start with the smoothed version
\[
F_X(s,\alpha) = \sum_{n=1}^\infty \frac{a(n)}{n^s} e^{-n^{1/d}z_X(\alpha)}
\]
of $F(s,\alpha)$, where $X>1$ is sufficiently large, $\alpha>0$ and
\begin{equation}
\label{3-1}
z_X(\alpha) = \frac{1}{X} + 2\pi i \alpha.
\end{equation}
Clearly, $F_X(s,\alpha)$ is absolutely convergent over $\CC$ and for every $\alpha>0$ we have
\[
\lim_{X\to\infty} F_X(s,\alpha) = F(s,\alpha) \qquad \text{for} \ \si>1.
\]
Our aim is to obtain a suitable expression for $F_X(s,\alpha)$ and then to investigate the limit as $X\to\infty$ for $s$ in certain regions inside the half-plane $\si<1$.

\smallskip
For a given $c>0$ and $-c<\si<2$, by Mellin's transform we have that 
\[
F_X(s,\alpha)  = \frac{1}{2\pi i} \int_{(d(c+2))} F(s+\frac{w}{d}) \Gamma(w) z_X(\alpha)^{-w} \d w.
\]
For $k\geq 0$ and $c$  sufficiently large in terms of $|s_k|$, we consider the ranges
\begin{equation}
\label{3-2}
\text{$\si\in\I_k$, where $\I_k$ is an arbitrary compact subinterval of $(-c,s_k)$}
\end{equation}
and, since $1/2 \leq s_0 \leq 1$, we define
\begin{equation}
\label{1-25}
k_0 = \min\{k\geq 1: s_k \leq 0\}.
\end{equation}
Moreover,  given $\delta>0$ sufficiently small in terms of $\I_k$ we write
\begin{equation}
\label{3-3}
u_k =
\begin{cases}
d(\delta-s_k) & \text{if} \ 0 \leq k \leq k_0-1 \\
\delta & \text{if} \ k \geq k_0
\end{cases}
\hskip-.1cm , \ \text{thus} \ \ \Re\big(s+\frac{u_k}{d}\big)<0 \ \ \text{for} \ \si\in\I_k.
\end{equation}
Since
\[
\text{$u_k<0$ for $0\leq k \leq k_0-1$ while $u_k>0$ for $k\geq k_0$,}
\]
for $s$ as in \eqref{3-2} we shift the line of integration in the above integral to $\Re(w) = u_k$, thus crossing the possible pole of $F(s+w/d)$ at $w= d(1-s)$ and, when $0\leq k \leq k_0-1$, also the simple poles of $\Gamma(w)$ at $w=-\nu$ with $0\leq \nu < ds_k$. Hence we get
\[
F_X(s,\alpha) = R_X(1-s,\alpha) + R_{k,X}(s,\alpha) +  \frac{1}{2\pi i} \int_{(u_k)} F(s+\frac{w}{d}) \Gamma(w) z_X(\alpha)^{-w} \d w,
\]
where $R_X(1-s,\alpha)$ is the residue of the integrand at $w= d(1-s)$ and, recalling that an empty sum equals 0, $R_{k,X}(s,\alpha)$ is the sums of the residues at $w=-\nu$ with $0\leq \nu < ds_k$. Next we apply the functional equation of $F(s)$ in the form \eqref{1-18}. In view of \eqref{3-3}, for $\si\in\I_k$ we may expand $ \overline{F}(1-s-w/d)$ and switch summation and integration, thus getting that
\begin{equation}
\label{3-4}
\begin{split}
F_X(s,\alpha) = \sum_{n=1}^\infty \frac{\overline{a(n)}}{n^{1-s}} &  \frac{1}{2\pi i} \int_{(u_k)} h_F \big(s+\frac{w}{d} \big) S_F \big(s+\frac{w}{d} \big) \Gamma(w) \big(\frac{z_X(\alpha)}{n^{1/d}}\big)^{-w} \d w \\
& +R_X(1-s,\alpha) + R_{k,X}(s,\alpha).
\end{split}
\end{equation}
For later use, we explicitly note that $R_{k,X}(s,\alpha)\equiv0$ if $k\geq k_0$.

\smallskip
In order to compute the residues we first note from \eqref{1-11} that
\[
 \sum_{m=1}^{m_F} \frac{d^m\gamma_m}{(w-d(1-s))^m}
\]
is the polar part of $F(s+w/d)$ at $w=d(1-s)$. Then obviously we have that
\begin{equation}
\label{3-5}
R_X(1-s,\alpha) = z_X(\alpha)^{-d(1-s)} \sum_{m=1}^{m_F} d^m\gamma_m \sum_{h=0}^{m-1} \frac{(-1)^h \log^h z_X(\alpha) \Gamma^{(m-1-h)}(d(1-s))}{h! (m-1-h)!}.
\end{equation}
Note that for any $\alpha>0$ and $s$ outside the singularities of $\Gamma^{(m-1-h)}(d(1-s))$ we have
\begin{equation}
\label{3-6}
\lim_{X\to\infty} R_X(1-s,\alpha) = R(1-s,\alpha),
\end{equation}
where $R(s,\alpha)$ is given by \eqref{1-19}. In particular, \eqref{3-6} holds for $\si<0$, and $R(1-s,\alpha)$ is meromorphic over $\CC$. We finally note that
\begin{equation}
\label{3-7}
R_{k,X}(s,\alpha) = \sum_{0\leq \nu < ds_k} \frac{(-1)^\nu}{\nu!} F\big(s-\frac{\nu}{d}\big) z_X(\alpha)^\nu,
\end{equation}
hence in particular $R_{k,X}(s,\alpha)$ is meromorphic over $\CC$.

\medskip
{\bf 3.2. Expansion of $h_F(s)$.}
As outlined in the Introduction, we reduce part of the treatment of the general case under consideration to the special case investigated in \cite{Ka-Pe/half}, which depends on the explicit computation of the Mellin-Barnes integral \eqref{1-6}. The first step in this direction is to obtain a suitable asymptotic expansion of type \eqref{1-22} for the function $h_F(s+w/d)$ in \eqref{3-4}. 

\smallskip
To this end, recalling \eqref{1-2}, we write
\begin{equation}
\label{3-8}
z=z(s,w) := \frac{d+1}{2} -ds -w = d(s_0-(s+\frac{w}{d})),
\end{equation}
so that
\[
s+\frac{w}{d} = \frac{d+1}{2d} -\frac{z}{d}
\]
and hence by \eqref{1-17}
\begin{equation}
\label{3-9}
\begin{split}
h_F \big(s+\frac{w}{d} \big) &=  \frac{\omega}{(2\pi)^r} Q^{1-2(s+w/d)} \prod_{j=1}^r \Gamma \big(\alpha_j+\frac{\lambda_j}{d} z \big) \Gamma \big(\beta_j+\frac{\lambda_j}{d} z \big) \\
&=  \frac{\omega}{(2\pi)^r} Q^{1-2(s+w/d)} \, \widetilde{h}_F \big(s+\frac{w}{d} \big),
\end{split}
\end{equation}
say, with
\begin{equation}
\label{3-10}
\alpha_j = \lambda_j \frac{d-1}{2d} +\overline{\mu}_j \quad \text{and} \quad \beta_j = 1- \lambda_j \frac{d+1}{2d} - \mu_j.
\end{equation}
Clearly, in view of \eqref{3-3} and \eqref{3-8}, for $s$ in the range \eqref{3-2} and $\Re(w)=u_k$ we have $\Re(z)>0$, thus by Stirling's formula we obtain
\[
\begin{split}
\log \widetilde{h}_F(s+\frac{w}{d}) &=  \sum_{j=1}^r \left( \log \Gamma(\alpha_j+\frac{\lambda_j}{d} z) + \log \Gamma(\beta_j+\frac{\lambda_j}{d} z) \right) \\
&= \sum_{j=1}^r (\alpha_j+\beta_j-1+\frac{2\lambda_j}{d}z) \log z + \sum_{j=1}^r (\alpha_j+\beta_j-1+\frac{2\lambda_j}{d}z) \log\frac{\lambda_j}{d} \\
&\hskip.5cm - z +r\log(2\pi) + \sum_{m=1}^M \frac{c_m}{z^m} + O(|z|^{-M-1})
\end{split}
\]
with any given integer $M\geq 0$ and certain constants $c_m\in\RR$. But, thanks to the hypothesis $\theta_F=0$, we have
\[
\sum_{j=1}^r (\alpha_j+\beta_j-1) = -\frac12,
\]
therefore
\[
 \sum_{j=1}^r (\alpha_j+\beta_j-1+\frac{2\lambda_j}{d}z) \log z = (z-\frac12) \log z
\]
and hence, again by Stirling's formula,
\[
\begin{split}
\log \widetilde{h}_F(s+\frac{w}{d}) &= \log \Gamma(z) + \sum_{j=1}^r (\alpha_j+\beta_j-1+\frac{2\lambda_j}{d}z) \log\frac{\lambda_j}{d} \\
&\hskip.5cm  + (r-\frac12)\log(2\pi) + \sum_{m=1}^M \frac{c'_m}{z^m} + O(|z|^{-M-1})
\end{split}
\]
with certain constants $c'_m\in\RR$. Moreover,
\[
\begin{split}
\sum_{j=1}^r (\alpha_j+\beta_j-1+\frac{2\lambda_j}{d}z) &\log\frac{\lambda_j}{d}  = \sum_{j=1}^r(-\frac{\lambda_j}{d} -2i\Im(\mu_j) +\frac{2\lambda_j}{d}z) \log\frac{\lambda_j}{d} \\
&= -\frac{1}{2d}\log\beta +\frac12\log d - 2i \sum_{j=1}^r \Im(\mu_j) \log \lambda_j + z\frac{\log\beta}{d} -z\log d,
\end{split}
\]
where
\begin{equation}
\label{3-11}
\beta=\prod_{j=1}^r \lambda_j^{2\lambda_j}.
\end{equation}

\smallskip
As a consequence, for $\Re(z)>0$ we have
\begin{equation}
\label{3-12}
\widetilde{h}_F(s+\frac{w}{d}) = (2\pi)^{r-1/2} \sqrt{d} \beta^{-1/(2d)} \prod_{j=1}^r \lambda_j^{-2i\Im(\mu_j)} \beta^{z/d} d^{-z} \Gamma(z) \Sigma_M(z),
\end{equation}
where
\begin{equation}
\label{3-13}
\begin{split}
\Sigma_M(z) :&= \exp\left( \sum_{m=1}^M \frac{c'_m}{z^m} +O(|z|^{-M-1})\right) \\
&=1 + \sum_{m=1}^M \frac{c''_m}{z^m} + O(|z|^{-M-1}) = \sum_{\ell=0}^M \frac{d_\ell}{(z-1)\cdots (z-\ell)} + E_M(z)
\end{split}
\end{equation}
with $c''_m\in\RR$, certain coefficients $d_\ell\in\RR$ and a function $E_M(z)$ satisfying
\begin{equation}
\label{3-14}
E_M(z)   = O(|z|^{-M-1}).
\end{equation}
Indeed, the last identity in \eqref{3-13} and bound \eqref{3-14} follow by elementary manipulations involving the geometric series, staring with $\ell=1$ and then recursively with $\ell=2,3 \dots$ Note that by convention the term with $\ell=0$ equals 1, thus 
\[
d_0=1 \ \ \text{when} \ \ \theta_F=0.
\]
Moreover, $E_M(z)$ is meromorphic for $\Re(z)>0$, and by \eqref{3-9},\eqref{3-12} and \eqref{3-13} its poles are at most at the points
\begin{equation}
\label{3-15}
z= \ell\in \{1,\dots,M\} \quad \text{and} \quad z= -\frac{d(\beta_j+h)}{\lambda_j}>0
\end{equation}
with $h\geq 0$ integer, $j=1,\dots,r$ and $\beta_j$ as in \eqref{3-10}. Note, by \eqref{3-13}, that the poles of the first type are at most simple. Therefore, recalling the definition of conductor $q$ and root number $\omega_F$ in \eqref{1-10}, from \eqref{3-8},\eqref{3-9},\eqref{3-11},\eqref{3-12},\eqref{3-13}, the factorial formula for the $\Gamma$ function and \eqref{1-2}, for $s$ in the range \eqref{3-2}, $\Re(w)=u_k$ and $M\geq0$ we have
\begin{equation}
\label{3-16}
h_F(s+\frac{w}{d}) =  \frac{\omega_F}{\sqrt{2\pi}} \left(\frac{q^{1/d}}{2\pi d}\right)^{\frac{d}{2}-d(s+\frac{w}{d})} \sum_{\ell=0}^M d_\ell \Gamma\big(d(s_\ell-(s+\frac{w}{d}))\big) + \widetilde{E}_M(s+\frac{w}{d}),
\end{equation}
say, $\widetilde{E}_M(s+w/d)$ being meromorphic for $\Re(z)>0$.

\smallskip
{\bf Remark 3.1.} Recall that we are working under the assumption $\theta_F=0$. Note that formally \eqref{3-16} corresponds to \eqref{1-22}, since $M\geq0$ is arbitrary. Moreover, in view of \eqref{3-12},\eqref{3-13} and \eqref{3-16} we have
\begin{equation}
\label{3-17}
\widetilde{E}_M(s+\frac{w}{d}) = \frac{\omega_F}{\sqrt{2\pi}} \left(\frac{q^{1/d}}{2\pi d}\right)^{\frac{d}{2}-d(s+\frac{w}{d})} \Gamma \big(d(s_0-(s+\frac{w}{d}))\big) E_M \big(d(s_0-(s+\frac{w}{d}))\big),
\end{equation}
thus, by \eqref{3-14}, $\widetilde{E}_M(s+w/d)$ is bounded by the modulus of the $M$-th term in \eqref{3-16} times $1/|s+w/d|$. Note also that the first set of potential poles of $E_M(z)$ in \eqref{3-15} is outside the range \eqref{3-2} when $\Re(w)=u_k$ and $\delta>0$ is sufficiently small, while in view of \eqref{1-15} the second set of such poles, if any, is cancelled by the zeros of $S_F(s+w/d)$, i.e.
\begin{equation}
\label{3-18}
E_M \big(d(s_0-(s+\frac{w}{d}))\big) S_F(s+\frac{w}{d})
\end{equation}
is holomorphic for $s$ and $\Re(w)$ as above. Note, however, that for $k\geq k_0$ the line $\si= s_k-\delta/d$ contains a potential pole of $E_M (d(s_0-(s+w/d))$ when $\Re(w)=u_k=\delta$. \qed

\medskip
{\bf 3.3. Computation of the integral in \eqref{3-4}.}
Now, recalling \eqref{1-15}, from \eqref{3-4}, \eqref{3-16} and \eqref{3-17} we have, for $s$ in the range \eqref{3-2}, that
\begin{equation}
\label{3-19}
\begin{split}
F_X(s,\alpha) &=  \frac{ \omega_F}{\sqrt{2\pi}} \left(\frac{q^{1/d}}{2\pi d}\right)^{\frac{d}{2}-ds} \sum_{\ell=0}^M d_\ell \sum_{j=-N}^N a_j e^{i\pi d\omega_j s} \sum_{n=1}^\infty \frac{\overline{a(n)}}{n^{1-s}}  \\
&\hskip1cm \times \frac{1}{2\pi i} \int_{(u_k)} \Gamma\big(d(s_\ell-s)-w\big) \Gamma(w) z_{j,X}(\alpha,n) ^{-w} \d w \\
& +  H_{M,X}^{(u_k)}(s,\alpha) + R_X(1-s,\alpha) + R_{k,X}(s,\alpha)
\end{split}
\end{equation}
where, recalling \eqref{3-1},
\begin{equation}
\label{3-20}
z_{j,X}(\alpha,n) = \frac{q^{1/d}z_X(\alpha) e^{-i\pi\omega_j}}{2\pi d n^{1/d}}
\end{equation}
and, after summation over $n$,
\begin{equation}
\label{3-21}
H_{M,X}^{(u_k)}(s,\alpha) = \frac{1}{2\pi i} \int_{(u_k)} \overline{F}\big(1-(s+\frac{w}{d})\big) S_F(s+\frac{w}{d}) \widetilde{E}_M\big(s+\frac{w}{d}\big)  \Gamma(w)  z_X(\alpha)^{-w}  \d w.
\end{equation}

\smallskip
Next we note that the integral in \eqref{3-19} is of the form \eqref{1-6}, but before applying formula \eqref{1-6} we have to ensure that conditions \eqref{1-7} are satisfied. In view of \eqref{3-20}, the second condition in \eqref{1-7} is always satisfied thanks to \eqref{1-16} since $\arg z_X(\alpha) = \pi/2 -\epsilon$ with some $\epsilon=\epsilon(X,\alpha)>0$. Moreover, in view of \eqref{3-3}, when $s$ is in the range \eqref{3-2} the first condition in \eqref{1-7} is satisfied for every $\ell=0,\dots,M$ provided
\[
k\geq k_0 \quad \text{and} \quad 0\leq M \leq k.
\]
When $0\leq k\leq k_0-1$ we have that $u_k<0$, so we shift the integral in \eqref{3-19} back to $\Re(w)=\delta$. For $0\leq \ell \leq k$ and $s$ in the range \eqref{3-2} we have $\Re(d(s_\ell-s)-w)>0$ whenever $u_k\leq \Re(w) \leq \delta$, hence we cross only the simple poles of $\Gamma(w)$ at $w=-\nu$ with $0\leq \nu < ds_k$. Thus we get
\begin{equation}
\label{3-22}
\begin{split}
\frac{1}{2\pi i} \int_{(u_k)} &\Gamma\big(d(s_\ell-s)-w\big) \Gamma(w) z_{j,X}(\alpha,n) ^{-w} \d w \\
&=  \frac{1}{2\pi i} \int_{(\delta)} \Gamma\big(d(s_\ell-s)-w\big) \Gamma(w) z_{j,X}(\alpha,n) ^{-w} \d w \\
&\hskip1cm  - \sum_{0\leq \nu < ds_k} \frac{(-1)^\nu}{\nu!} \Gamma\big(d(s_\ell-s)+\nu\big) z_{j,X}(\alpha,n) ^\nu \\
& = \frac{1}{2\pi i} \int_{(\delta)} \Gamma\big(d(s_\ell-s)-w\big) \Gamma(w) z_{j,X}(\alpha,n)^{-w} \d w - \widetilde{R}_{k,\ell,j,X}(s,\alpha,n),
\end{split}
\end{equation}
say. Therefore, choosing
\[
M=k,
\]
recalling \eqref{3-3} and observing that $\widetilde{R}_{k,\ell,j,X}(s,\alpha,n) \equiv0$ if $k\geq k_0$, we may apply formula \eqref{1-6} to the integrals over the line $\Re(w)=\delta$ in \eqref{3-19} and \eqref{3-22} for every $k\geq0$, $0\leq \ell\leq k$ and $s$ in the range \eqref{3-2}, thus obtaining that
\begin{equation}
\label{3-23}
\begin{split}
\frac{1}{2\pi i} \int_{(u_k)} & \Gamma\big(d(s_\ell-s)-w\big)  \Gamma(w) z_{j,X}(\alpha,n) ^{-w} \d w \\
&= \Gamma\big(d(s_\ell-s)\big) \big(1+ z_{j,X}(\alpha,n)\big)^{d(s-s_\ell)} - \widetilde{R}_{k,\ell,j,X}(s,\alpha,n).
\end{split}
\end{equation}

\smallskip
Suppose now that $\si\in\I_k\cap (-\infty,-2\delta)$, where $\I_k$ is defined in \eqref{3-2}. Then, substituting \eqref{3-23} into \eqref{3-19} with $M=k$, the sum over $n$ is absolutely convergent for evey $k\geq 0$. Hence we get
\begin{equation}
\label{3-24}
\begin{split}
F_X(s,\alpha) =  &\frac{\omega_F}{\sqrt{2\pi}} \left(\frac{q^{1/d}}{2\pi d}\right)^{\frac{d}{2}-ds} \sum_{\ell=0}^k d_\ell \Gamma\big(d(s_\ell -s)\big)  \overline{F}^*_{\ell,X}(1-s,\alpha)  \\
& +R_X(1-s,\alpha)  + R_{k,X}(s,\alpha) - \widetilde{R}_{k,X}(s,\alpha) +  H_{k,X}^{(u_k)}(s,\alpha),
\end{split}
\end{equation}
where $R_X(1-s,\alpha)$, $R_{k,X}(s,\alpha)$ and $H_{k,X}^{(u_k)}(s,\alpha)$ are given by \eqref{3-5}, \eqref{3-7} and \eqref{3-21},
\begin{equation}
\label{3-25}
\overline{F}^*_{\ell,X}(1-s,\alpha) = \sum_{j=-N}^N a_j e^{i\pi d\omega_j s} \sum_{n=1}^\infty \frac{\overline{a(n)}}{n^{1-s}} (1+z_{j,X}(\alpha,n))^{d(s-s_\ell)}
\end{equation}
and
\begin{equation}
\label{3-26}
\begin{split}
\widetilde{R}_{k,X}(s,\alpha) = e^{as+b} &  \sum_{0\leq \nu < ds_k} \frac{(-1)^\nu}{\nu!}  \sum_{\ell=0}^k d_\ell  \Gamma\big(d(s_\ell-s)+\nu\big) \\
&\times \sum_{j=-N}^N a_j e^{i\pi d\omega_j s} \sum_{n=1}^\infty \frac{\overline{a(n)}}{n^{1-s}}  z_{j,X}(\alpha,n) ^\nu
\end{split}
\end{equation}
with 
\begin{equation}
\label{3-27}
e^{as+b}= \frac{\omega_F}{\sqrt{2\pi}} \left(\frac{q^{1/d}}{2\pi d}\right)^{\frac{d}{2}-ds}
\end{equation}
and $z_{j,X}(\alpha,n)$, $\omega_j$, $a_j$ given by \eqref{3-20} and \eqref{1-15}. Note that the argument of $1+z_{j,X}(\alpha,n)$ lies in $[-\epsilon, \pi-\epsilon]$, with $\epsilon\to0$ as $X\to\infty$. But a computation shows that \eqref{3-26} transforms to
\[
\begin{split}
\widetilde{R}_{k,X}(s,\alpha) = & \sum_{0\leq \nu<ds_k}  \frac{(-1)^\nu}{\nu!}  \overline{F}\big(1-(s-\frac{\nu}{d})\big) S_F\big(s-\frac{\nu}{d}\big) \\
&\times \left\{ \frac{ \omega_F}{\sqrt{2\pi}}  \left(\frac{q^{1/d}}{2\pi d}\right)^{\frac{d}{2}-d(s-\frac{\nu}{d})}  \sum_{\ell=0}^k d_\ell  \Gamma\Big(d\big(s_\ell- (s-\frac{\nu}{d})\big)\Big) \right\} z_X(\alpha)^\nu.
\end{split}
\]
Hence, thanks to the functional equation \eqref{1-18}, \eqref{3-7} and \eqref{3-16} with $w=-\nu$ and $M=k$ we get
\begin{equation}
\label{3-28}
\begin{split}
\widetilde{R}_{k,X}(s,\alpha)  =  R_{k,X}(s,\alpha) - E_{k,X}(s,\alpha),
\end{split}
\end{equation}
where
\begin{equation}
\label{3-29}
E_{k,X}(s,\alpha) =  \sum_{0\leq \nu<ds_k}  \frac{(-1)^\nu}{\nu!}  \overline{F}\big(1-(s-\frac{\nu}{d})\big) S_F\big(s-\frac{\nu}{d}\big) \widetilde{E}_k\big(s-\frac{\nu}{d}\big) z_X(\alpha)^\nu.
\end{equation}
Note that $E_{k,X}(s,\alpha) \equiv0$ if $k\geq k_0$, and that $E_{k,X}(s,\alpha)$ is well defined for $\si\in\I_k$ in view of \eqref{3-16} and \eqref{3-17}, since $\Re(z)=\Re(d(s_0-s)+\nu)>0$ in such a range. 

\smallskip
Writing
\begin{equation}
\label{3-30}
H_{k,X}(s,\alpha) = 
\begin{cases}
H_{k,X}^{(u_k)}(s,\alpha) + E_{k,X}(s,\alpha) & \text{if} \ 0\leq k \leq k_0-1 \\
H_{k,X}^{(\delta)}(s,\alpha) & \text{if} \ k\geq k_0,
\end{cases}
\end{equation}
from \eqref{3-3}, \eqref{3-24}, \eqref{3-28} and \eqref{3-30} we finally have for $k\geq 0$ and $\si\in\I_k\cap (-\infty,-2\delta)$ that
\begin{equation}
\label{3-31}
\begin{split}
F_X(s,\alpha) =  &\frac{\omega_F}{\sqrt{2\pi}} \left(\frac{q^{1/d}}{2\pi d}\right)^{\frac{d}{2}-ds} \sum_{\ell=0}^k d_\ell \Gamma\big(d(s_\ell -s)\big)  \overline{F}^*_{\ell,X}(1-s,\alpha)  \\
& +R_X(1-s,\alpha) +  H_{k,X}(s,\alpha).
\end{split}
\end{equation}
Note that the only terms in \eqref{3-31} requiring that $\si$ belongs to the subinterval $\I_k\cap (-\infty,-2\delta)$, instead of the whole interval $\I_k$, are the functions $\overline{F}^*_{\ell,X}(1-s,\alpha)$. However, these terms are close to the functions $\overline{F}_\ell(s,\alpha)$, which are known to be entire by Theorem 1, so this issue will not be a problem in the next section.

\smallskip
{\bf Remark 3.2.} Again, recall that we are working under the assumption that $\theta_F=0$, although such an assumption is not relevant for this remark. For $0\leq k \leq k_0-1$ we could shift the line of integration in $H_{k,X}^{(u_k)}(s,\alpha)$ back to $\Re(w)=\delta$, to obtain that
\[
H_{k,X}^{(u_k)}(s,\alpha) = - E_{k,X}(s,\alpha) + H_{k,X}^{(\delta)}(s,\alpha).
\]
Indeed, for $\si\in\I_k\cap (-\infty,-2\delta)$ and $u_k \leq \Re(w) \leq \delta$ the possible pole of $\overline{F}(1-s-w/d)$ is not crossed, and $\widetilde{E}_k(s+w/d)$ is well defined for $\si\in\I_k$ since  $\Re(z)=\Re(d(s_0-s)-w)>0$. Hence \eqref{3-31} would assume the uniform shape
\[
\begin{split}
F_X(s,\alpha) =  &\frac{\omega_F}{\sqrt{2\pi}} \left(\frac{q^{1/d}}{2\pi d}\right)^{\frac{d}{2}-ds} \sum_{\ell=0}^k d_\ell \Gamma\big(d(s_\ell -s)\big)  \overline{F}^*_{\ell,X}(1-s,\alpha)  \\
& +R_X(1-s,\alpha) +  H_{k,X}^{(\delta)}(s,\alpha)
\end{split}
\]
for every $k\geq 0$. If $k=k_0-1$, $s_{k_0}<0$ and $s_{k_0} < \si<0$, then the argument developed in the next section still allows to let $X\to\infty$ in $H_{k_0-1,X}^{(\delta)}(s,\alpha)$. However, in the next section we consider the intervals $s_{k+1}< \si<s_k$ for every $k\geq0$. But if $0\leq k\leq k_0-1$ and $\si>0$, then letting $X\to\infty$ in $H_{k,X}^{(\delta)}(s,\alpha)$ involves some difficulties due to the poor known bounds on the Lindel\"of $\mu$-function of $F(s)$ in the critical strip; see the computations in Lemma 3.1 below. Such a limit could instead be performed if $F(s)$ satisfies the Lindel\"of Hypothesis. \qed

\medskip
{\bf 3.4. Limit as $X\to\infty$.}
We first remark, as in Section 2.2 of \cite{Ka-Pe/half}, that letting $X\to\infty$ in \eqref{3-31} requires some care. Indeed, as we already pointed out, the limit of $F_X(s,\alpha)$ is $F(s,\alpha)$ when $\si>1$, but \eqref{3-31} holds in the range $\I_k\cap (-\infty,-2\delta)$. Moreover, the limit of the terms $(1+ z_{j,X}(\alpha,n))^{d(s-s_\ell)}$ in \eqref{3-25} is not always well defined, since $1+ z_{j,X}(\alpha,n)$ may vanish as $X\to\infty$. 

\smallskip
As in Section 2.2 of \cite{Ka-Pe/half}, we first compute $F_X(s,\alpha)$ in a different way. Since $F_X(s,\alpha)$ is the twist of $F(s,\alpha)$ by $e^{-n^{1/d}/X}$, by Mellin's transform we have that
\[
F_X(s,\alpha) = \frac{1}{2\pi i} \int_{(c')} F(s+\frac{w}{d},\alpha) \Gamma(w) X^w \d w
\]
for $s$ as in \eqref{3-2}, where $c'=c'(k)>0$ is sufficiently large. Then we shift the line of integration to $\Re(w) = -\delta$, thus crossing the simple pole at $w=0$ with residue $F(s,\alpha)$ coming from $\Gamma(w)$, and the possible simple poles at $w=d(s_\ell - s)$ with $0\leq \ell\leq k$ coming from $F(s+w/d,\alpha)$, with residues
\begin{equation}
\label{3-32}
\kappa_{\ell,X}(s,\alpha) = d \rho_\ell(\alpha) \Gamma(d(s_\ell-s)) X^{d(s_\ell-s)}.
\end{equation}
Hence for $s$ in the range \eqref{3-2} we have
\begin{equation}
\label{3-33}
\begin{split}
F_X(s,\alpha) &= F(s,\alpha) +  \sum_{\ell=0}^k \kappa_{\ell,X}(s,\alpha)+ \frac{1}{2\pi i} \int_{(-\delta)} F(s+\frac{w}{d},\alpha) \Gamma(w) X^w \d w \\
&= F(s,\alpha) + \Sigma_X(s,\alpha) + I_X(s,\alpha),
\end{split}
\end{equation}
say.
Note that $\Sigma_X(s,\alpha)$ vanishes if $\alpha\not\in$ Spec$(F)$, since $F(s,\alpha)$ is entire in this case. Moreover, it is immediate to see that
\begin{equation}
\label{3-34}
\lim_{X\to\infty} I_X(s,\alpha) = 0
\end{equation}
for every $\alpha>0$ and $s$ as in \eqref{3-2}.

\smallskip
Now, recalling \eqref{1-13}, \eqref{1-16} and \eqref{3-20}, we observe that
\begin{equation}
\label{3-35}
1+ z_{-N,X}(\alpha,n_\alpha) = \frac{iq^{1/d}}{2\pi dn_\alpha^{1/d}} \frac{1}{X}
\end{equation}
and
\begin{equation}
\label{3-36}
\lim_{X\to\infty} (1+ z_{j,X}(\alpha,n)) = 
\begin{cases}
1+ e^{i\pi(\frac12-\omega_j)} \left(\frac{n_\alpha}{n}\right)^{1/d} \neq 0 \quad &\text{if} \ j \neq -N \ \ \text{or} \ \ n\neq n_\alpha, \\
0 \quad &\text{if} \ j =- N \ \ \text{and} \ \ n=n_\alpha.
\end{cases}
\end{equation}
Hence in view of \eqref{3-35} and \eqref{3-36} we rewrite \eqref{3-25} as
\begin{equation}
\label{3-37}
\begin{split}
\overline{F}^*_{\ell,X}(1-s,\alpha) &= \sum_{j=-N}^N a_j e^{i\pi d\omega_j s} \sideset{}{^\flat}\sum_{n\geq1} \frac{\overline{a(n)}}{n^{1-s}} (1+z_{j,X}(\alpha,n))^{d(s-s_\ell)} \\
&\hskip1.5cm+ a_{-N} e^{-i\frac{\pi}{2} ds} \frac{\overline{a(n_\alpha)}}{n_\alpha^{1-s}} \Big(\frac{iq^{1/d}}{2\pi dn_\alpha^{1/d}} \frac{1}{X}\Big)^{d(s-s_\ell)} \\
&= \overline{F}_{\ell,X}(1-s,\alpha) + a_{-N} e^{-i\frac{\pi}{2} ds} \frac{\overline{a(n_\alpha)}}{n_\alpha^{1-s}} \Big(\frac{iq^{1/d}}{2\pi dn_\alpha^{1/d}} \frac{1}{X}\Big)^{d(s-s_\ell)},
\end{split}
\end{equation}
say, where the symbol $^\flat$ has the same meaning as in \eqref{1-20}. Note, recalling \eqref{1-12}, that if $\alpha \not\in$ Spec$(F)$ then the last term in the right hand side of \eqref{3-37} vanishes and we may omit the symbol $^\flat$ in the second sum. Therefore, for $\si\in\I_k\cap (-\infty,-2\delta)$ and $k\geq0$, from \eqref{3-31}, \eqref{3-33} and \eqref{3-37} we have
\begin{equation}
\label{3-38}
\begin{split}
F&(s,\alpha) + \Sigma_X(s,\alpha) + I_X(s,\alpha) = R_X(1-s,\alpha) + H_{k,X}(s,\alpha) \\
& +  \frac{\omega_F}{\sqrt{2\pi}} \left(\frac{q^{1/d}}{2\pi d}\right)^{\frac{d}{2}-ds} \sum_{\ell=0}^k d_\ell \Gamma(d(s_\ell-s)) \overline{F}_{\ell,X}(1-s,\alpha) + \widetilde{\Sigma}_X(s,\alpha),
\end{split}
\end{equation}
where, recalling \eqref{3-27},
\begin{equation}
\label{3-39}
\widetilde{\Sigma}_X(s,\alpha)  =  a_{-N} e^{as+b} e^{-i\frac{\pi}{2} ds} \frac{\overline{a(n_\alpha)}}{n_\alpha^{1-s}} \sum_{\ell=0}^k d_\ell \Gamma(d(s_\ell -s)) \left(\frac{iq^{1/d}}{2\pi dn_\alpha^{1/d}} \frac{1}{X}\right)^{d(s-s_\ell)}.
\end{equation}
Moreover, the terms $\Sigma_X(s,\alpha)$ and $\widetilde{\Sigma}_X(s,\alpha)$ vanish unless $\alpha\in$ Spec$(F)$.

\smallskip
Now we are ready to let $X\to\infty$. We first deal with the functions $\overline{F}_{\ell,X}(1-s,\alpha)$, proceeding as in the proof of Theorem 1; so we only outline the argument. We start splitting $\overline{F}_{\ell,X}(1-s,\alpha)$ in two parts as in \eqref{2-1} and choose 
\[
\rho = d(s-s_\ell) \quad \text{and} \quad z= z_{j,X}(n,\alpha),
\]
then we proceed to obtain the analog of \eqref{2-3} with a fixed $R$ ($R=2([d]+1)$ is already sufficient). Since all the involved sums are now either finite or absolutely convergent for $\si\leq s_0$, we may let $X\to\infty$ and then, in view of \eqref{1-20}, \eqref{2-2} and \eqref{3-36}, we proceed backwards to reconstruct $\overline{F}_\ell(1-s,\alpha)$. In this way we obtain that for $\si$ in the range \eqref{3-2} and $0\leq \ell\leq k$
\begin{equation}
\label{3-40}
\lim_{X\to \infty} \overline{F}_{\ell,X}(1-s,\alpha) = \overline{F}_\ell(1-s,\alpha).
\end{equation}
Moreover, we already noticed, see \eqref{3-6} and \eqref{3-34}, that 
\begin{equation}
\label{3-41}
\lim_{X\to \infty} R_X(1-s,\alpha) = R(1-s,\alpha) \quad \text{and} \quad \lim_{X\to \infty} I_X(s,\alpha) = 0
\end{equation}
for every $\alpha>0$ and $s$ as in \eqref{3-2}. Note that $R(1-s,\alpha)$ is holomorphic in the range \eqref{3-2}.

\smallskip
Next we deal with $H_{k,X}(s,\alpha)$ and, in view of \eqref{3-17}, \eqref{3-21}, \eqref{3-27}, \eqref{3-29} and \eqref{3-30}, we define
\begin{equation}
\label{3-42}
\begin{split}
H_k^{(u_k)}(s,\alpha) &=  e^{as+b}  \frac{1}{2\pi i} \int_{(u_k)} \overline{F}\big(1-s-\frac{w}{d}\big) \Gamma\big(d(s_0-s)-w\big) \\
&\times E_k\big(d(s_0-s)-w\big) S_F\big(s+\frac{w}{d}\big)  \Gamma(w) \left(\frac{q^{1/d} \alpha i}{d}\right)^{-w} \d w,
\end{split}
\end{equation}
\begin{equation}
\label{3-43}
E_{k}(s,\alpha) =  \sum_{0\leq \nu<ds_k}  \frac{(-1)^\nu}{\nu!}  \overline{F}\big(1-(s-\frac{\nu}{d})\big) S_F\big(s-\frac{\nu}{d}\big) \widetilde{E}_k\big(s-\frac{\nu}{d}\big) (2\pi i\alpha)^\nu
\end{equation}
and
\begin{equation}
\label{3-44}
H_{k}(s,\alpha) = 
\begin{cases}
H_{k}^{(u_k)}(s,\alpha) + E_{k}(s,\alpha) & \text{if} \ 0\leq k \leq k_0-1 \\
H_{k}^{(\delta)}(s,\alpha) & \text{if} \ k\geq k_0.
\end{cases}
\end{equation}
Note that \eqref{3-44} formally corresponds to $H_{k,\infty}(s,\alpha)$, see \eqref{3-30}.

\medskip
{\bf Lemma 3.1.} {\sl Let $F(s)$ be as in Theorem 2 with $\theta_F=0$, and let $k\geq 0$. Then $H_k(s,\alpha)$ is holomorphic in the vertical strip $s_{k+1}<\si<s_k$ and}
\[
\lim_{X\to\infty} H_{k,X}(s,\alpha) = H_k(s,\alpha).
\]

\medskip
{\it Proof.} Let $I_k(s,X)$ and $I_k(s)$ denote the integrals in \eqref{3-21} and \eqref{3-42}, respectively. We shall prove that $I_k(s,X)$ and $I_k(s)$ converge absolutely and uniformly for $s$ in any compact subset $\K$ of the strip $(s_{k+1},s_k)$, the convergence of $I_k(s,X)$ being uniform also with respect to $X$. In particular, this allows to switch limit and integration, thus showing that
\[
\lim_{X\to\infty} H_{k,X}(s,\alpha) = H_k(s,\alpha),
\]
since clearly $E_{k,X}(s,\alpha)$ tends to $E_k(s,\alpha)$. Moreover, this also shows that $H_k(s,\alpha)$ is holomorphic for $s_{k+1}<\si<s_k$. Indeed, the integrand in \eqref{3-42} is holomorphic at least for $s_{k+1}<\si<s_k-\delta/d$ with any sufficiently small $\delta>0$ thanks to Remark 3.1 with $M=k$, see \eqref{3-18}, and the same holds for $E_k(s,\alpha)$ for similar reasons. In view of the sharp similarity of such integrals, we treat them simultaneously.

\smallskip
Let $\K$ be as above, $s\in\K$, the value of $\delta>0$ in \eqref{3-3} be sufficiently small in terms of $\K$, and let $\I$ be the projection of $\K$ on the real axis. We write  $w=u_k+iv$ and split both $I_k(s,X)$ and $I_k(s)$ in two parts $I_1$ and $I_2$ with $|v|\leq 1$ and $|v|>1$, respectively. We have
\[
z_X(\alpha)^{-w}, i^{-w} \ll e^{\frac{\pi}{2}|v|}
\]
uniformly in $X$ and, in view of the definition of $u_k$ in \eqref{3-3},
\[
\overline{F}(1-s-\frac{w}{d}) \ll_\delta 1, \qquad S_F(s+\frac{w}{d}) \ll e^{\frac{\pi}{2}|dt +v|}
\]
uniformly for $\si\in\I$. We also have
\[
\Gamma(w) \ll \frac{1}{\delta+|v|} \qquad \text{for} \ \ |v|\leq 1
\]
and, by \eqref{3-14} with $M=k$,
\[
 E_k(d(s_0-s)-w) \ll_\delta (1+|dt+v|)^{-k-1},
\]
again uniformly for $\si\in\I$. Then, by Stirling's formula, in view of \eqref{1-2} we get
\begin{equation}
\label{3-45}
I_1 \ll_\delta \int_{|v|\leq 1} (1+ |dt+v|)^{d(s_0-\si) - u_k - \frac12 - k -1} \frac{\d v}{\delta+|v|}  \ll_{\K} 1
\end{equation}
uniformly in $X$ and for $s\in\K$. Moreover, for a certain $v_0=v_0(\K)$ we have
\begin{equation}
\label{3-46}
\begin{split}
I_2 &\ll_\delta \int_{|v|>1} (1+|dt+v|)^{d(s_0-\si) - u_k - \frac12 - k -1} (1+|v|)^{u_k-\frac12} \d v \\
&\ll_{\K} \int_{|v|\geq v_0} (1+|v|)^{d(s_{k+1}-\si)+1} \d v \ll_{\K} 1
\end{split}
\end{equation}
uniformly in $X$ and for $s\in\K$, since $d(s_{k+1}-\si)>0$. The lemma follows now from \eqref{3-45} and \eqref{3-46}. \qed

\medskip
Finally we deal with the remaining two terms in \eqref{3-38}, namely $\Sigma_X(s,\alpha)$ and $\widetilde{\Sigma}_X(s,\alpha)$. Recall that here we may assume that $\alpha\in$ Spec$(F)$, otherwise both terms vanish. Thanks to \eqref{3-33} and \eqref{3-39} we rewrite such quantities as
\begin{equation}
\label{3-47}
\Sigma_X(s,\alpha) = \sum_{\ell=0}^k  a_\ell(s,\alpha) X^{d(s_\ell-s)} \quad \text{and} \quad \widetilde{\Sigma}_X(s,\alpha) = \sum_{\ell=0}^k  \widetilde{a}_\ell(s,\alpha) X^{d(s_\ell-s)},
\end{equation}
where $a_\ell(s,\alpha)$ and $  \widetilde{a}_\ell(s,\alpha)$ are independent of $X$. Since we already proved that all  the other terms in \eqref{3-38} tend to finite limits as $X\to\infty$, we have that
\[
\Sigma_X(s,\alpha) = \widetilde{\Sigma}_X(s,\alpha) + O(1) \qquad \text{as} \ \ X\to\infty;
\]
here the error term $O(1)$ depends on $s$ and other parameters, but is bounded in $X$. Moreover, for $s$ in the range \eqref{3-2} the exponent of $X$ in \eqref{3-47} satisfies $\Re(d(s_\ell-s))>0$, thus from \eqref{3-38} and \eqref{3-47} we deduce, starting with $\ell=0$ and then recursively $\ell=1,2,\dots$, that 
\begin{equation}
\label{3-48}
\Sigma_X(s,\alpha) = \widetilde{\Sigma}_X(s,\alpha).
\end{equation}

\smallskip
Now we can conclude the proof of Theorem 2. Indeed, from \eqref{3-40}, \eqref{3-41}, Lemma 3.1 and \eqref{3-48} we may let $X\to\infty$ in \eqref{3-38} to obtain that for $s_{k+1}<\si<s_k$
\begin{equation}
\label{3-49}
F(s,\alpha) =  \frac{\omega_F}{\sqrt{2\pi}} \left(\frac{q^{1/d}}{2\pi d}\right)^{\frac{d}{2}-ds} \sum_{\ell=0}^k d_\ell \Gamma\big(d(s_\ell-s) \big) \overline{F}_{\ell}(1-s,\alpha) + R(1-s,\alpha) + H_k(s,\alpha),
\end{equation}
all the involved functions being holomorphic in such range. Theorem 2, when $\theta_F=0$, follows now from \eqref{3-49}, \eqref{1-2}, \eqref{1-24} and the bounds for $H_k(s,\alpha)$ we are going to prove in the next section, see Lemma 3.2.

\smallskip
{\bf Remark 3.3.} Note that the uniform convergence argument in Lemma 3.1 still holds for $H^{(\delta)}_{k,X}(s,\alpha)$ and $H^{(\delta)}_k(s,\alpha)$ if 
\[
k=k_0-1, \quad  s_{k_0}<0 \quad \text{and} \quad s_{k_0} < \si<-\eta,
\] 
$\eta>0$ arbitrary. Hence by Remark 3.2 and the above limit process we have that \eqref{3-49} holds with $H_k(s,\alpha)$ replaced by $H^{(\delta)}_k(s,\alpha)$ in the above cases. This remark will be used in the next section. \qed

\medskip
{\bf 3.5. Bounds for $H_k(s,\alpha)$.}
The required bounds for $H_k(s,\alpha)$ are obtained by a refinement of Lemma 3.1, coupled with the Phragm\'en-Lindel\"of theorem and \eqref{3-49}. Once more, we recall that $k_0\geq1$ is given by \eqref{1-25}, thus $s_{k_0}\leq 0$, $-1< ds_{k_0} \leq 0$ and $s_{k_0-1}>0$.

\medskip
{\bf Lemma 3.2.} {\sl Let $F(s)$ be as in Theorem 2 with $\theta_F=0$, and let $k\geq 0$. Then $H_k(s,\alpha)$ is meromorphic over $\CC$, all its poles lie in a horizontal strip of bounded height and has polynomial growth on every vertical strip. Moreover, let $\ep>0$ be arbitrarily small and 
\[
m_\ep=\max(d(s_{k_0}+\ep),-1/2).
\] 
Then as $|t|\to\infty$ we have}
\[
H_k(s,\alpha) \ll_{\ep,k}
\begin{cases}
|t|^{-1/2+\ep} & \text{if} \ k\geq k_0+1, \ s_{k+1}\leq \si\leq s_k; \\
|t|^{-1/2+\ep} & \text{if} \ k=k_0, \ s_{k_0}=0, \ s_{k_0+1}\leq \si\leq -\ep; \\
|t|^{m_\ep+\ep} & \text{if} \ k_0-1\leq k\leq k_0, \ s_{k_0}<0, \ \si\in(s_{k_0}-\ep,-\ep]\cap[s_{k+1},s_k].
\end{cases}
\]

\medskip
{\it Proof.} Let $k\geq0$. It is an immediate consequence of \eqref{3-49}, the properties of $F(s,\alpha)$, Theorem 1 and \eqref{1-19} that $H_k(s,\alpha)$ is meromorphic over $\CC$ with all poles in a horizontal strip of bounded height and polynomial growth on every vertical strip.

\smallskip
Now, given a small $\eta>0$, we first assume that 
\begin{equation}
\label{3-50}
k\geq k_0 \quad  \text{and} \quad  s_{k+1} + \frac{1}{2d} \leq \si\leq s_k-\eta,
\end{equation}
and let $I_k(s)$, $I_1$ and $I_2$ be as in Lemma 3.1. We choose $\delta=\eta/2$ in \eqref{3-3}, thus  $\Re(1-s-w/d)>1$ and hence
\[
\overline{F}\big(1-s-\frac{w}{d}\big) \ll_\eta 1.
\]
Moreover, recalling Remark 3.1, by \eqref{3-14} we have
\[
E_k\big(d(s_0-s)-w\big) S_F\big(s+\frac{w}{d}\big) \ll_{\eta,k}  (1+|dt+v|)^{-k-1}  e^{\frac{\pi}{2}|dt +v|}
\]
and by Stirling's formula
\[
\Gamma\big(d(s_0-s)-w\big) \ll_k  (1+|dt+v|)^{d(s_0-\si)-\delta-1/2}  e^{-\frac{\pi}{2}|dt +v|}.
\]
Hence, similarly as in Lemma 3.1, the bound \eqref{3-45} becomes
\[
I_1 \ll_{\eta,k} \int_{|v|\leq 1} (1+ |dt+v|)^{d(s_{k+1}-\si)-\delta-1/2} \frac{\d v}{\delta+|v|} \ll_{\eta,k} (1+|t|)^{-1}
\]
for $\si$ as in \eqref{3-50}, while \eqref{3-46} becomes, with the notation of Lemma B in the Appendix,
\[
I_2 \ll_{\eta,k}  \int_{|v|>1} (1+ |dt+v|)^{d(s_{k+1}-\si)-\delta-1/2} (1+|v|)^{\delta-1/2} \d v \ll_{\eta,k} I(dt,A,B)
\]
with
\[
A= d(\si-s_{k+1}) + \frac12 +\delta \quad \text{and} \quad B = \frac12 -\delta.
\]
Thus from Lemma B we have
\[
I_2 \ll_{\eta,k} (1+|t|)^{-1/2+\eta}
\]
for $\si$ as in \eqref{3-50}, and therefore as $|t|\to\infty$
\begin{equation}
\label{3-51}
H_k(s,\alpha) \ll_{\eta,k} |t|^{-1/2+\eta}
\end{equation}
under conditons \eqref{3-50}.

\smallskip
In order to deal with the interval
\[
s_{k+1} \leq \si < s_{k+1} + \frac{1}{2d}
\]
we observe that, with the notation \eqref{3-27}, \eqref{3-49} gives for any $k\geq0$
\begin{equation}
\label{3-52}
H_k(s,\alpha) = H_{k+1}(s,\alpha) + e^{as+b} d_{k+1} \Gamma\big(d(s_{k+1}-s) \big) \overline{F}_{k+1}(1-s,\alpha).
\end{equation}
Thus from \eqref{3-51} with $\eta$ replaced by $\eta/d$ and applied to $H_{k+1}(s,\alpha)$ on the line $\si=s_{k+1}-\eta/d$, Theorem 1 and Stirling's formula, thanks to \eqref{3-52} we get
\begin{equation}
\label{3-53}
H_k(s_{k+1}-\eta/d+it,\alpha) \ll_{\eta,k} |t|^{-1/2+\eta}
\end{equation}
as $|t|\to\infty$. Since $H_k(s,\alpha)$ has polynomial growth on every vertical strip and all its poles lie in a horizontal strip of bounded height, from \eqref{3-51} applied on the line $\si= s_{k+1}+1/(2d)$, \eqref{3-53} and a standard argument based on the Phragm\'en-Lindel\"of theorem we obtain that
\begin{equation}
\label{3-54}
H_k(s,\alpha) \ll_{\eta,k} |t|^{-1/2+\eta} \quad \text{for} \quad s_{k+1}\leq \si \leq s_k-\eta
\end{equation}
as $|t|\to\infty$, provided $k\geq k_0$. 

\smallskip
Next we deal with the interval 
\[
s_k-\eta \leq \si \leq s_k.
\]
If $k\geq 1$ we rewrite \eqref{3-52} in the form
\begin{equation}
\label{3-55}
H_k(s,\alpha) = H_{k-1}(s,\alpha) - e^{as+b} d_{k} \Gamma\big(d(s_{k}-s) \big) \overline{F}_{k}(1-s,\alpha)
\end{equation}
and, for $k\geq k_0+1$, we apply \eqref{3-54} to $H_k(s,\alpha)$ on the line $\si=s_k-\eta$ and to $H_{k-1}(s,\alpha)$ on the line $\si=s_k$. Thus, as before, by \eqref{3-55}, Theorem 1, Stirling's formula and the Phragm\'en-Lindel\"of theorem we get as $|t|\to\infty$ that
\begin{equation}
\label{3-56}
H_k(s,\alpha) \ll_{\eta,k} |t|^{-1/2+\eta} \quad \text{for} \quad s_{k+1}\leq \si \leq s_k
\end{equation}
provided $k\geq k_0+1$. Therefore, \eqref{3-54} and \eqref{3-56} conclude the proof of the lemma if $s_{k_0}=0$.

\smallskip
Suppose finally that $s_{k_0}<0$ and write $k^*=k_0-1$. In this case, thanks to Remark 3.3, in the range $s_{k_0}<\si\leq-\eta$ we may use $H_{k^*}^{(\delta)}(s,\alpha)$ in place of $H_{k^*}(s,\alpha)$, and the same argument leading to \eqref{3-51} gives
\begin{equation}
\label{3-57}
H_{k^*}^{(\delta)}(-\eta+it,\alpha) \ll_{\eta} |t|^{\max(d(s_{k_0}+\eta),-1/2)+\eta}
\end{equation}
as $|t|\to\infty$. Then, thanks to \eqref{3-52} with $k=k^*$ and \eqref{3-54} applied with $k=k_0$, $\eta$ replaced by $\eta/d$ and $\si=s_{k_0}-\eta/d$, the above argument based on Theorem 1 and Stirling's formula gives
\begin{equation}
\label{3-58}
H_{k^*}^{(\delta)}(s_{k_0}-\eta/d+it,\alpha) \ll_{\eta} |t|^{-1/2+\eta}
\end{equation}
as $|t|\to\infty$. Hence \eqref{3-57}, \eqref{3-58} and the argument based on the Phragm\'en-Lindel\"of theorem give, as $|t|\to\infty$, that
\begin{equation}
\label{3-59}
H_{k^*}^{(\delta)}(s,\alpha) \ll_{\eta} |t|^{\max(d(s_{k_0}+\eta),-1/2)+\eta} \quad \text{for} \quad s_{k_0}\leq \si\leq -\eta.
\end{equation}
Thus the same bound holds for $H_{k^*}(s,\alpha)$ thanks to Remark 3.3. Moreover, a similar argument shows that the bound in \eqref{3-59} holds also for $H_{k_0}(s,\alpha)$ in the range $s_{k_0}-\eta<\si\leq s_{k_0}$. Therefore the lemma follows from \eqref{3-54}, \eqref{3-56} and \eqref{3-59} if $s_{k_0}<0$. \qed

\smallskip
{\bf Remark 3.4.} We first note that, clearly, the constants in the symbol $\ll_{\eta,k}$ in Lemma 3.2 depend also on $\alpha$ and $F(s)$. One can get a slightly sharper bound for $H_k(s,\alpha)$ by refining the arguments in Lemma 3.2, where $|t|^\ep$ is replaced by a power of $\log|t|$. Essentially, the point is to choose $\eta$ and $\delta$ in the proof as suitable functions of $|t|$, to use the full force of Lemma B in the Appendix and to employ a variant of the Phragm\'en-Lindel\"of theorem. Moreover, the exponent in the third bound in Lemma 3.2 may possibly be improved to $-1/2+\ep$. However, such improvements do not lead to any improvements in our present applications. We finally note that bounds for $H_k(s,\alpha)$ in the range $-\ep < \si \leq s_0$ can also be obtained, but in this case the estimates involve the Lindel\"of $\mu$-function of $F(s)$. \qed

\medskip
{\bf 3.6. Dropping the assumption $\theta_F=0$.}
Finally we drop the assumption $\theta_F=0$, and for simplicity we write 
\begin{equation}
\label{3-60}
s^*= s-i\theta_F.
\end{equation}
We first observe that if $F(s)$ satisfies functional equation \eqref{1-9} with data $\omega, Q, \lambda_j,\mu_j$, then the shifted function 
\[
G(s) = F(s^*) = \sum_{n=1}^\infty \frac{a^*(n)}{n^s}, \hskip1.5cm a^*(n) = a(n)n^{i\theta_F},
\]
although not formally a member of $\S^\sharp$ if $F(s)$ has a pole at $s=1$, satisfies \eqref{1-9} with data
\begin{equation}
\label{3-61}
\omega^*=\omega Q^{2i\theta_F}, \quad Q^*=Q, \quad \lambda_j^* = \lambda_j, \quad \mu_j^* = \mu_j - i\theta_F\lambda_j.
\end{equation}
Hence degree and conductor of $G(s)$ remain unchanged, while in view of \eqref{1-10} we have
\begin{equation}
\label{3-62}
\xi_G = \xi_F -id\theta_F = \eta_F \quad \text{and} \quad \omega_G = \omega_F \big(Q^2 \prod_{j=1}^r\lambda_j^{2\lambda_j}\big)^{i\theta_F} = \omega_F d^{id\theta_F} \big(\frac{q^{1/d}}{2\pi d}\big)^{id\theta_F}.
\end{equation}
Thus $G(s)$ has $\theta_G=0$, therefore we can apply the arguments in Sections  3.1 to 3.5 to $G(s)$ and its standard twist
\[
G(s,\alpha) = F(s^*,\alpha) = \sum_{n=1}^\infty \frac{a^*(n)}{n^s} e(-\alpha n^{1/d}).
\]
Note, by \eqref{1-14}, that the potential poles of $G(s)$ and $G(s,\alpha)$ are now at the points
\[
s= 1+i\theta_F \qquad \text{and} \qquad s=s_\ell \  \text{with} \   \ell=0,1,2,\dots
\]
with polar parts
\[
\sum_{m=1}^{m_F} \frac{\gamma_m}{(s-(1+i\theta_F))^m} \qquad \text{and} \qquad \frac{\rho_\ell(\alpha)}{(s-s_\ell)},
\]
respectively. In what follows we briefly sketch the changes to be made with respect to the case $\theta_F=0$. 

\smallskip
As at the beginning of Section 3.1 we express $G_X(s,\alpha)$ by a Mellin trasform and then shift the integration line to $\Re(w) = u_k$. Since now we cross the potential pole of $G(s+w/d)$ at $w=d(1-s+i\theta)=d(1-s^*)$, see \eqref{3-60} for $s^*$, equation \eqref{3-4} becomes
\[
\begin{split}
G_X(s,\alpha) = \sum_{n=1}^\infty \frac{\overline{a^*(n)}}{n^{1-s}} &  \frac{1}{2\pi i} \int_{(u_k)} h_G \big(s+\frac{w}{d} \big) S_G \big(s+\frac{w}{d} \big) \Gamma(w) \big(\frac{z_X(\alpha)}{n^{1/d}}\big)^{-w} \d w \\
& + R_X(1-s^*,\alpha) + R_{k,X}(s,\alpha),
\end{split}
\]
where $R_X(s,\alpha)$ is as in \eqref{1-19} and $R_{k,X}(s,\alpha)$ is defined by \eqref{3-7} in terms of $G(s)$. Moreover, by \eqref{3-61} or recalling that $G(s)=F(s^*)$, expression \eqref{1-15} for $G(s)$ is
\begin{equation}
\label{3-63}
S_G(s) = \sum_{j=-N}^N a_j e^{i\pi d \omega_js^*} = \sum_{j=-N}^N a_j e^{\pi d \omega_j \theta_F} e^{i\pi d \omega_js} =  \sum_{j=-N}^N a_j^* e^{i\pi d \omega_js},
\end{equation}
say, and in view of \eqref{3-62} equation \eqref{3-16} in Section 3.2 is replaced by
\[
\begin{split}
h_G(s+\frac{w}{d}) =  \frac{\omega_F}{\sqrt{2\pi}} \left(\frac{q^{1/d}}{2\pi d}\right)^{\frac{d}{2}-d(s^*+\frac{w}{d})} \sum_{\ell=0}^M d_\ell \Gamma\big(d(s_\ell-(s+\frac{w}{d}))\big)  +  \widetilde{E}_M(s+\frac{w}{d}).
\end{split}
\]
Clearly, here $\widetilde{E}_M(s+w/d)$ is defined by \eqref{3-17} in terms of $G(s)$, and 
\[
d_\ell = d^{id\theta_F} d_\ell(G)
\] 
where the $d_\ell(G)$ come from the argument in Section 3.2 applied to $G(s)$. In particular, in the general case we have
\begin{equation}
\label{3-64}
d_0 = d^{id\theta_F}.
\end{equation}
As a consequence, \eqref{3-31} in Section 3.3 is now replaced by
\begin{equation}
\label{3-65}
G_X(s,\alpha) =  e^{as^*+b} \sum_{\ell=0}^k d_\ell \Gamma\big(d(s_\ell -s)\big)  \overline{G}^*_{\ell,X}(1-s,\alpha)  +R_X(1-s^*,\alpha) +  H_{k,X}(s,\alpha),
\end{equation}
where $H_{k,X}(s,\alpha)$ is defined by \eqref{3-30} in terms of $G(s)$, $a$ and $b$ are as in \eqref{3-27} and, in view of \eqref{3-63}, $\overline{G}^*_{\ell,X}(1-s,\alpha)$ is defined by \eqref{3-25} with $a_j$ and $a(n)$ replaced by $a_j^*$ and $a^*(n)$, respectively.

\smallskip
The limit as $X\to\infty$ in \eqref{3-65} is performed as in Section 3.4. First, recalling that the potential poles of $G(s,\alpha)$ are at the points $s=s_\ell$, we have that  equation \eqref{3-38} changes to
\begin{equation}
\label{3-66}
\begin{split}
G(s&,\alpha) + \Sigma_X(s,\alpha) + I_X(s,\alpha) = R_X(1-s^*,\alpha) +   H_{k,X}(s,\alpha) \\
& +  e^{as^*+b}  \sum_{\ell=0}^k d_\ell \Gamma(d(s_\ell-s)) \overline{G}_{\ell,X}(1-s,\alpha) +  \widetilde{\Sigma}_X(s,\alpha).
\end{split}
\end{equation}
Again, here $I_X(s,\alpha)$ is now defined in terms of $G(s)$, $ \overline{G}_{\ell,X}(1-s,\alpha)$ is defined by \eqref{3-37} with $a^*_j$ and $a^*(n)$ in place of $a_j$ and $a(n)$, \eqref{3-39} is replaced by
\begin{equation}
\label{3-67}
\widetilde{\Sigma}_X(s,\alpha)  =  a_{-N}^* e^{as^*+b} e^{-i\frac{\pi}{2} ds} \frac{\overline{a^*(n_\alpha)}}{n_\alpha^{1-s}} \sum_{\ell=0}^k d_\ell \Gamma(d(s_\ell -s)) \left(\frac{iq^{1/d}}{2\pi dn_\alpha^{1/d}} \frac{1}{X}\right)^{d(s-s_\ell)}
\end{equation}
and $\Sigma_X(s,\alpha)$ is unchanged, hence by \eqref{3-32}
\begin{equation}
\label{3-68}
 \Sigma_X(s,\alpha) =\sum_{\ell=0}^k  d \rho_\ell(\alpha) \Gamma(d(s_\ell-s)) X^{d(s_\ell-s)}.
\end{equation}
Moreover, as in \eqref{3-48} we now have
\begin{equation}
\label{3-69}
 \Sigma_X(s,\alpha) =\widetilde{\Sigma}_X(s,\alpha).
\end{equation}
Then, letting $X\to\infty$ in \eqref{3-66} we finally obtain
\begin{equation}
\label{3-70}
G(s,\alpha) =  e^{as^*+b}  \sum_{\ell=0}^k d_\ell \Gamma(d(s_\ell-s)) \overline{G}_{\ell}(1-s,\alpha) + R(1-s^*,\alpha) + H_{k}(s,\alpha),
\end{equation}
where once again $H_k(s,\alpha)$ is defined by \eqref{3-44} in terms of $G(s)$, and
\begin{equation}
\label{3-71}
\overline{G}_\ell(1-s,\alpha) = \sum_{j=-N}^N a_j^* e^{i\pi d\omega_j s} \sideset{}{^\flat} \sum_{n\geq1} \frac{\overline{a^*(n)}}{n^{1-s}} \left(1+ e^{i\pi(\frac12-\omega_j)} \left(\frac{n_\alpha}{n}\right)^{1/d} \right)^{d(s-s_\ell)}.
\end{equation}

\smallskip
Now we shift $s$ to $s+i\theta_F$, thus $s^*$ in \eqref{3-60} goes back to $s$. From \eqref{1-14}, \eqref{1-20}, the definition of $a^*(n)$, \eqref{3-63} and \eqref{3-71} we have that \eqref{3-70} becomes
\begin{equation}
\label{3-72}
F(s,\alpha) =  e^{as+b}  \sum_{\ell=0}^k d_\ell \Gamma(d(s_\ell^*-s)) \overline{F}_{\ell}(1-s,\alpha) + R(1-s,\alpha) + H_{k}(s+i\theta_F,\alpha).
\end{equation}
Theorem 2 follows  in the general case from \eqref{3-72} by renaming $H_k(s,\alpha)$ the function $H_{k}(s+i\theta_F,\alpha)$, which is not given explicitly in Theorem 2 and whose properties are clearly the same as those stated in Lemmas 3.1 and 3.2.

\medskip
{\bf 3.7. Proof of Theorem 3.}
Let $\alpha\in$ Spec$(F)$. Writing $\widetilde{\Sigma}_X(s,\alpha)$ and $\Sigma_X(s,\alpha)$ in \eqref{3-67} and \eqref{3-68} in the form \eqref{3-47}, from \eqref{3-69} we deduce that
\[
a_\ell(s,\alpha) = \widetilde{a}_\ell(s,\alpha), \hskip1.5cm \ell = 0,\dots, k.
\]
Moreover, from \eqref{3-68} we have
\[
a_\ell(s,\alpha) = d \rho_\ell(\alpha) \Gamma(d(s_\ell-s))
\]
and from \eqref{3-27} and \eqref{3-67}
\[
 \widetilde{a}_\ell(s,\alpha) = a_{-N}^* \frac{\omega_F}{\sqrt{2\pi}}  \left(\frac{q^{1/d}}{2\pi d}\right)^{\frac{d}{2}-ds+id\theta_F} e^{-i\frac{\pi}{2} ds} \frac{\overline{a^*(n_\alpha)}}{n_\alpha^{1-s}} d_\ell \Gamma(d(s_\ell -s)) \left(\frac{iq^{1/d}}{2\pi dn_\alpha^{1/d}} \right)^{d(s-s_\ell)}.
\]
But a computation based on \eqref{1-15}, \eqref{3-62} and \eqref{3-63} shows that 
\[
a_{-N}^* =  e^{-i\frac{\pi}{2} \xi_G} = e^{-i\frac{\pi}{2} \eta_F},
\]
hence by \eqref{1-10} and \eqref{1-14} we get
\[
\begin{split}
\rho_\ell(\alpha) &= \frac{d_\ell}{d} \frac{\omega_F e^{-i\frac{\pi}{2} \eta_F}}{\sqrt{2\pi}}  \left(\frac{q^{1/d}}{2\pi d}\right)^{\frac{d}{2}-ds+id\theta_F} e^{-i\frac{\pi}{2} ds} \frac{\overline{a(n_\alpha)}}{n_\alpha^{1-s+i\theta_F}} \left(\frac{iq^{1/d}}{2\pi dn_\alpha^{1/d}} \right)^{d(s-s_\ell)} \\
&= \frac{d_\ell}{d} \frac{\omega_F}{\sqrt{2\pi}} e^{-i\frac{\pi}{2}(\xi_F+ds_\ell^*)} \left(\frac{q^{1/d}}{2\pi d}\right)^{\frac{d}{2}-ds_\ell^*} \frac{\overline{a(n_\alpha)}}{n_\alpha^{1-s_\ell^*}},
\end{split}
\]
and Theorem 3 follows.

\bigskip
\section{Proof of Theorem 4 and some remarks}

\smallskip
In the first two sections we prove Theorem 4, the third one is devoted to several remarks on the strict functional equation. Recall that we assume $\theta_F=0$.

\medskip
{\bf 4.1. (i) $\Longleftrightarrow$ (ii).}
Assuming (i), we have that (ii) follows at once from \eqref{1-29} and \eqref{4-19} below. 

\smallskip
Next, assuming that (ii) holds, we show that \eqref{1-25} holds with $H_k(s,\alpha)\equiv0$ for every $k\geq h$ and all $\alpha>0$. Then, since $d_h\neq0$ and $\overline{F}_h(s,\alpha)\not\equiv0$ by Theorem 1, $h$ must be minimal as required by Definition 1.1. Moreover, since by Theorem 3 we have $d_\ell=0$ for every $\ell \geq h+1$, if $H_h(s,\alpha)\equiv0$ for every $\alpha>0$ then $H_k(s,\alpha)\equiv0$ for every $k\geq h$ and $\alpha>0$. Hence it remains to prove that $H_h(s,\alpha)\equiv0$ for every $\alpha>0$.

\smallskip
From \eqref{3-13},\eqref{3-14} and the argument leading to \eqref{3-16} and \eqref{3-17}, with $w=0$ and $z=d(s_0-s)$, we have that
\begin{equation}
\label{4-1}
h_F(s)  = e^{as+b} \Gamma(z) \big(P(z) + E_h(z)\big).
\end{equation}
Here $e^{as+b}$ is as in \eqref{3-27}, $P(z)$ is defined by
\[
P(z) = 1 + \sum_{\ell=1}^h \frac{d_\ell}{(z-1)\cdots (z-\ell)}
\]
and we recall that $E_h(z)$ is meromorphic for $\Re(z)>0$, see \eqref{3-15}, and satisfies
\begin{equation}
\label{4-2}
E_h(z) = O\big(|z|^{-A}\big)
\end{equation}
for every $A>0$ in view of \eqref{3-14}, since $d_\ell=0$ for every $\ell\geq h+1$. Our aim is therefore to show that $E_h(z)\equiv0$, which immediately implies (i) thanks to \eqref{3-42},\eqref{3-43},\eqref{3-17} and the above observations. To this end we compute the asymptotic expansion of $h_F(s)$ in Section 3.2 using Lemma D in the Appendix instead of Stirling's formula. We set $w=0$ and follow the steps from \eqref{3-8} to \eqref{3-12}, pointing out only the changes needed in the error terms due to the application of Lemma D, since the treatment of the main terms is exactly the same.

\smallskip
We apply Lemma D to the logarithm of each $\Gamma$-factor in \eqref{3-9} with $\delta>0$ arbitrarily small,
\[
 \rho= \frac{\lambda_jz}{d}, \quad a= \alpha_j \ \text{or} \  a=\beta_j \quad \text{and} \ \ \gamma = \frac{d}{\lambda_j}\eta \ \ \text{with} \  \eta=\min_{j=1,\dots,r} \frac{\lambda_j}{d}, \ \text{say};
 \]
note that $0<\gamma\leq1$.  Thus for $|\arg(z)+\beta|<\pi(1/2-\delta)$ we obtain an error of the form
 \[
 \int_0^{\frac{d}{\lambda_j}\eta e^{i\beta}} \psi(w,a) e^{-w\lambda_jz/d} \frac{\d w}{w} + O_\delta(e^{-\kappa |z|})
 \]
with some $\kappa>0$ and $\beta$ arbitrarily close to $\pm\pi/2$. Then, by the change of variable $\frac{\lambda_j}{d} w \mapsto w$, this error becomes
 \[
 \int_0^{\eta e^{i\beta}} \psi(\frac{d}{\lambda_j}w,a) e^{-wz} \frac{\d w}{w} + O_\delta(e^{-\kappa |z|})
 \]
 again with some $\kappa>0$. Therefore, after summation over $j=1,\dots,r$ and the other main term computations in Section 3.2 leading to \eqref{3-12}-\eqref{3-16}, we obtain that
\begin{equation}
\label{4-3}
\log\Big(\frac{e^{-as-b}}{\Gamma(z)} h_F(s)\Big) = \int_0^{\eta e^{i\beta}} \psi_F(w) e^{-wz} \frac{\d w}{w} + O_\delta(e^{-\kappa |z|})
\end{equation}
for $|\arg(z)+\beta|<\pi(1/2-\delta)$, where
\[
\psi_F(w) = \sum_{m=1}^\infty c(m) w^m \hskip1.5cm \text{for $|w|<2\pi\eta$}
\]
with certain coefficients $c(m)$ depending on $F(s)$. Hence from \eqref{4-3} and Lemma C in the Appendix, with $P(z)$ as above, we deduce that
\begin{equation}
\label{4-4}
\log\Big(\frac{e^{-as-b}}{\Gamma(z)P(z)} h_F(s)\Big) = \int_0^{\eta e^{i\beta}} \big(\psi_F(w) - \psi_P(w) \big) e^{-wz} \frac{\d w}{w} + O_\delta(e^{-\kappa |z|}),
\end{equation}
where $\psi_P(w)$ is defined by \eqref{6-4}, and therefore for  $|w|<2\pi\eta$
\[
\psi_F(w) - \psi_P(w) = \sum_{m=1}^\infty e(m) w^m
\] 
with certain coefficients $e(m)$.

\smallskip
Suppose now that there exists $m\geq 1$ with $e(m)\neq0$, and let $m_0$ be the least such $m$. Then, thanks to the above conditions on $\arg(z)$ and $\beta$, we have
\[
\begin{split}
\int_0^{\eta e^{i\beta}} \big(\psi_F(w) - \psi_P(w) \big) e^{-wz} \frac{\d w}{w} &= e(m_0) \int_0^{\eta e^{i\beta}} e^{-wz} w^{m_0-1} \d w \\
&\hskip.5cm + O\Big( \sum_{m>m_0} |e(m)| \int_0^{\eta e^{i\beta}} e^{-|zw|\delta} |w|^{m-1} |\d w|\Big) \\
&= \frac{e(m_0)}{z^{m_0}} \int_0^{z\eta e^{i\beta}} e^{-\xi} \xi^{m_0-1} \d \xi + O_{\delta}\Big( \frac{1}{|z|^{m_0+1}} \Big) \\
& = \frac{e(m_0) \Gamma(m_0)}{z^{m_0}} + + O_{\delta}\Big( \frac{1}{|z|^{m_0+1}} \Big).
\end{split}
\]
As a consequence, \eqref{4-4} becomes
\[
\log\Big(\frac{e^{-as-b}}{\Gamma(z)P(z)} h_F(s)\Big) = \frac{c_0}{z^{m_0}} + O_{\delta}\Big( \frac{1}{|z|^{m_0+1}} \Big)
\]
with a certain constant $c_0\neq0$, thus
\begin{equation}
\label{4-5}
\begin{split}
h_F(s) &= s^{as+b} \Gamma(z) P(z) \exp\Big(  \frac{c_0}{z^{m_0}} + O_{\delta}\Big( \frac{1}{|z|^{m_0+1}} \Big) \Big) \\
& = s^{as+b} \Gamma(z) P(z) \Big(1 +  \frac{c_0}{z^{m_0}} + O_{\delta}\Big( \frac{1}{|z|^{m_0+1}} \Big) \Big).
\end{split}
\end{equation}
Comparing \eqref{4-1},\eqref{4-2} and \eqref{4-5} we deduce that $c_0=0$, a contradiction unless
\[
\psi_F(w) = \psi_P(w). 
\]
Therefore  \eqref{4-4} becomes
\begin{equation}
\label{4-6}
\log\Big(\frac{e^{-as-b}}{\Gamma(z)P(z)} h_F(s)\Big) =  O_{\delta'}(e^{-\kappa |z|})
\end{equation}
for $|z|$ sufficiently large, $|\arg(z)|<\pi-\delta'$ and $\delta'>0$ arbitrarily small.

\smallskip
{\bf Remark 4.1.} Note that the conditions involving $\beta$ are not anymore present here. Indeed, the whole range $|\arg(z)|<\pi-\delta'$ is covered by means of two overlapping angular regions of the form $|\arg(z)+\beta|<\pi(1/2 - \delta)$ with $\beta$ arbitrarily close to $\pm\pi/2$, on each of which \eqref{4-6} holds. \qed

\smallskip
Now we can conclude the proof by the arguments on p.106-107 of \cite{Ka-Pe/2002a}. Let $z=c_0+it$ with a given large $c_0>0$ and let
\[
f(t) = \frac{e^{-as-b}h_F(s)}{\Gamma(z)P(z)} -1.
\]
Thanks to \eqref{4-6}, for $k=0,1,\dots$ we have
\begin{equation}
\label{4-7}
\mu_k(f) = \int_{-\infty}^{+\infty} f(t) t^k \d t = \frac{1}{i^{k+1}} \int_{(c_0)} \Big(\frac{e^{-as-b}h_F(s)}{\Gamma(z)P(z)} - 1\Big) (z-c_0)^k \d z =0,
\end{equation}
as it can be checked by shifting the line of integration to $+\infty$. Moreover, since again by \eqref{4-6} we have $f(t) \ll e^{-\kappa|t|}$, the Fourier transform
\[
\widehat{f}(z) = \int_{-\infty}^{+\infty} f(\xi) e^{-\xi z} \d \xi
\]
of $f(t)$ is holomorphic for $|\Im(z)|<\kappa$. But
\[
\widehat{f}^{(k)}(0) = (-2\pi i)^k \mu_k(f) = 0
\]
for every $k\geq 0$ by \eqref{4-7}, hence $\widehat{f}(z)\equiv 0$. Thus $f(t)\equiv0$ as well, and therefore
\[
h_F(s)  = e^{as+b} \Gamma(z) P(z),
\]
i.e. $E_h(z)\equiv 0$ by \eqref{4-1} and (i) follows. \qed

\medskip
{\bf 4.2. (i) $\Longleftrightarrow$ (iii).}
Assume (i). Then by \eqref{1-28} we know that $d_\ell=0$ for $\ell\geq h+1$ and hence, with the notation in \eqref{3-27}, thanks to \eqref{1-2},\eqref{1-14},\eqref{1-22},\eqref{1-23} and \eqref{4-19} we have
\begin{equation}
\label{4-8}
h_F(s) = e^{as+b} \sum_{\ell=0}^h d_\ell \Gamma\big(\frac{d+1}{2} - \ell - ds \big),  \hskip1.5cm d_0d_h\neq0.
\end{equation}
Thus by \eqref{1-17} we have
\begin{equation}
\label{4-9}
\prod_{j=1}^r \big(\Gamma(\lambda_j(1-s)+\overline{\mu}_j) \Gamma(1-\lambda_js-\mu_j)\big) = e^{a's+b'} \sum_{\ell=0}^h d_\ell \Gamma\big(\frac{d+1}{2} - \ell - ds \big)
\end{equation}
with certain $a'\in\RR$ and $b'\in\CC$. The poles of the left hand side coincide with the generalized arithmetic progressions
\begin{equation}
\label{4-10}
\frac{1-\mu_j+\nu}{\lambda_j} \quad \text{and} \quad 1 + \frac{\overline{\mu}_j+\nu}{\lambda_j} \qquad \text{with $j=1,\dots,r$ and $\nu=0,1,\dots$,}
\end{equation}
with steps $1/\lambda_j$, while those of the right hand side are contained in the progression
\begin{equation}
\label{4-11}
\frac{d+1}{2d}  + \frac{\nu}{d} \qquad \text{with $\nu\geq -h$,}
\end{equation}
with steps $1/d$. Hence we deduce that
\[
\text{$\mu_j\in\RR$ for every $j$} \quad \text{and} \quad \lambda_j = \frac{d}{\nu_j} \ \text{with some $\nu_j\in\NN$, for every $j$,}
\] 
thus in particular
\[
\frac{\lambda_i}{\lambda_j} \in\QQ \quad \text{for every $i,j$.}
\]
Therefore, applying the multiplication formula of the $\Gamma$ function if necessary, we may suppose without loss of generality that all $\lambda_j$ are equal, and hence there exits an integer $N\geq1$ such that $F(s)$ has a $\gamma$-factor of the form
\begin{equation}
\label{4-12}
\gamma(s) = Q^s \prod_{j=1}^N \Gamma\Big(\frac{d}{2N} s +\mu_j\Big).
\end{equation}
Assuming now that the $\lambda_j$ are as in \eqref{4-12}, taking $\nu=0$ in \eqref{4-10} and comparing with \eqref{4-11} we deduce that there exist integers $n_j$ such that
\begin{equation}
\label{4-13}
\mu_j = \frac{1}{2N} \big(n_j - \frac{d+1}{2}\big)= \frac{2n_j-d-1}{4N} \hskip1.5cm j=1,\dots,r
\end{equation}
and $n_j\geq (d+1)/2$, since $\Re(\mu_j)\geq 0$. Moreover, the left hand side of  \eqref{4-9} becomes
\begin{equation}
\label{4-14}
\prod_{j=1}^N \Big( \Gamma \Big(\frac{d}{2N}(1-s) + \frac{n_j}{2N} - \frac{d+1}{4N}\Big) \Gamma \Big(\frac{d}{2N}(1-s) + \frac{2N+1-n_j}{2N} - \frac{d+1}{4N}\Big)\Big)
\end{equation}
and hence all its poles must be simple in accordance with the right hand side. Clearly, this restriction holds if and only if conditions (a) and (b) in \eqref{1-31} are satisfied. Together with \eqref{4-12} and \eqref{4-13}, this proves the first part of (iii), and it remains to show that $N_0=h$.

\smallskip
Inserting the values of the $n_j$ coming from the right hand side of \eqref{1-32} into \eqref{4-14}, rearranging terms if necessary and using the factorial formula for the $\Gamma$ function, \eqref{4-14} becomes
\begin{equation}
\label{4-15}
Q(z)\prod_{j=0}^{N-1} \Big( \Gamma \big(z+\frac{n_0+j}{2N}\big) \Gamma \big(z + \frac{n_0+j}{2N} + \frac12\big)\Big),
\end{equation}
where $z=d(1-s)/(2N) - (d+1)/(4N)$ and $Q(z)\not\equiv0$ is a certain real monic polynomial of degree $N_0\geq0$. Next we apply the duplication formula to \eqref{4-15}, thus obtaining from \eqref{4-9} that
\[
h_F(s) = e^{a''s+b''} Q(z) \prod_{j=0}^{N-1} \Gamma\Big(2z+ \frac{n_0+j}{N}\Big)
\]
with certain $a''\in\RR$ and $b''\in\CC$. Hence an application of the multiplication formula gives
\begin{equation}
\label{4-16}
h_F(s) = e^{As+B} Q(z) \Gamma(2Nz + n_0)
\end{equation}
with certain $A\in\RR$ and $B\in\CC$. Now we write $Q(z)$ in the form
\[
Q(z) = b_0 + b_1w + b_2w(w+1) + \cdots + b_{N_0} w(w+1)\cdots(w+N_0-1)
\]
with $b_j\in\RR$ and $w = 2Nz+n_0$. Hence, applying the factorial formula again, \eqref{4-16} takes the form
\[
\begin{split}
h_F(s) &= e^{As+B} \sum_{\ell=0}^{N_0} b_\ell \Gamma\big(d(1-s) +n_0+\ell - \frac{d+1}{2}\big) \\
&= e^{As+B} \sum_{\ell=0}^{N_0} b_\ell \Gamma\big(\frac{d+1}{2} +n_0-1+\ell - ds\big) \\
&= e^{As+B} \sum_{\ell=1-n_0-N_0}^{1-n_0}b_{1-n_0-\ell} \Gamma\big(\frac{d+1}{2} -\ell - ds\big).
\end{split}
\]
But, summing over $j$ the integers inside both sets in \eqref{1-32} and equating the results, one easily checks that $1-n_0=N_0$. Hence we finally get
\begin{equation}
\label{4-17}
h_F(s) = e^{As+B} \sum_{\ell=0}^{N_0}b_{N_0-\ell} \Gamma\big(\frac{d+1}{2} -\ell - ds\big),
\end{equation}
and comparing \eqref{4-17} with \eqref{4-8} we deduce that $N_0=h$, thus (iii) follows.

\smallskip
Suppose now that (iii) holds. Then (i) follows reversing the above arguments. Indeed, we start with the definition \eqref{1-17} of $h_F(s)$ and substitute the values of $\lambda_j$ and $\mu_j$ coming from (iii), thus getting that the $\Gamma$-factors in \eqref{1-17} transform to \eqref{4-14}. Then we insert in \eqref{4-14} the values of the $n_j$ coming from the right hand side of \eqref{1-32}, and follow the same computations leading to \eqref{4-15},\eqref{4-16} and finally to \eqref{4-17} with certain $A'$ and $B'$ in place of $A$ and $B$. But we know that $h_F(s)$ has the form \eqref{1-22}, hence $h_F(s)$ must have the shape \eqref{4-8} with
\[
h=N_0, \quad a=A', \quad  e^bd_\ell =e^{B'}b_{N_0-\ell}.
\]
This shows that (i) holds. Indeed, \eqref{4-8} implies that $H_h(s,\alpha)\equiv0$ for every $\alpha>0$, as it can be seen from the arguments in Section 3.2 and definition \eqref{3-44} since $E_h(z)\equiv0$ in this case. Choosing $\alpha\in$ Spec$(F)$ we deduce that the poles of $F(s,\alpha)$ are at most at the points $s_\ell$ with $0\leq \ell\leq h$, thus by Theorem 3 we have $d_\ell=0$ for $\ell\geq h+1$. But this implies that $H_k(s,\alpha)\equiv0$ for every $k\geq h$ and $\alpha>0$, as we have already seen at the beginning of Section 4.1. Finally, $h$ is minimal since we have chosen $N_0$ as the minimal integer such that \eqref{1-32} holds. \qed

\medskip
{\bf 4.3. Remarks on the strict functional equation.} 
Here we assume that $F(s)$ is as in Theorem 2, and for simplicity that $\theta_F=0$.

\smallskip
{\bf 1.} Suppose that $F(s,\alpha)$ satisfies a strict functional equation. Then, comparing \eqref{1-26} and \eqref{1-30}, by meromorphic continuation we have that if $h\geq 1$ then for $0\leq k\leq h-1$, $\alpha>0$ and $s\in\CC$
\begin{equation}
\label{4-18}
H_k(s,\alpha) =  \frac{\omega_F}{\sqrt{2\pi}} \left(\frac{q^{1/d}}{2\pi d}\right)^{\frac{d}{2}(1-2s)} \sum_{\ell = k+1}^{h}d_{\ell} \Gamma\big(d(s_{\ell}-s) \big) \overline{F}_{\ell}(1-s,\alpha).
\end{equation}
Since $h$ in Definition 1.1 is minimal, if $F(s,\alpha)$ satisfies (i) of Theorem 4 we also have that
\begin{equation}
\label{4-19}
d_h\neq 0,
\end{equation}
hence thanks to \eqref{1-28} we can write
\[
h = \max \{\ell\geq 0: d_\ell\neq0\}.
\]

\smallskip
{\bf 2.} Suppose that $H_k(s,\alpha_0)\equiv0$ for $k=h$ and some $\alpha_0\in$ Spec$(F)$. Then one can prove that $H_k(s,\alpha_0)\equiv 0$ for every $k\geq h$. This shows that the conditions in the definition of strict functional equation may be somewhat relaxed. However, the strong form in Definition 1.1 is preferable in view of the equivalences in Theorem 4.

\smallskip
{\bf 3.} {\it If $F(s)$ has a pole at $s=1$ and $F(s,\alpha)$ satisfies a strict functional equation, then} $m_F=1$, i.e. the pole of $F(s)$ is simple. Indeed, by \eqref{1-9} the polar order $m_F$ is at most the order of pole of $\gamma(s)$ at $s=0$. But the $\gamma$-factors in (iii) of Theorem 4 have only simple poles, hence our assertion follows. This follows also directly from \eqref{1-30}, since $F(s,\alpha)$ and $\Gamma\big(d(1-s)-e_\ell\big)$ have at most simple poles and the functions $\overline{F}_\ell(1-s,\alpha)$ are entire by Theorem 1. Hence $R(1-s,\alpha)$ has at most simple poles as well, thus $m_F\leq 1$ in view of \eqref{1-19}.

\smallskip
{\bf 4.} {\it If $F(s)$ has a pole at $s=1$ and $F(s,\alpha)$ satisfies a strict functional equation, then the degree $d$ is an odd integer}. In particular, the degree conjecture is true in this case. Indeed, for $\alpha\not\in$ Spec$(F)$ the left hand side of \eqref{1-30} is entire, hence the poles of $R(1-s,\alpha)$ have to cancel with the poles of $\Gamma\big(d(1-s)-e_\ell\big) \overline{F}_\ell(1-s,\alpha)$. Hence, recalling \eqref{1-19} and \eqref{1-24}, we must have
\[
\frac{d-1}{2} \in \NN,
\]
and our assertion follows.

\smallskip
{\bf 5.} Suppose that $F(s)$ has a pole at $s=1$ and $F(s,\alpha)$ satisfies a strict functional equation, and let $\alpha\not\in$ Spec$(F)$ and
\[
\rho_F= \res_{s=1} F(s).
\]
Hence $m_F=1$ and $d\geq1$ is an odd integer by the previous remarks, thus $e_\ell\in\NN$ by \eqref{1-24}, and the left hand side of \eqref{1-30} is entire. Therefore the poles of $R(1-s,\alpha)$ at
\[
s=1+\frac{\nu}{d} \qquad \text{with} \  \nu\geq0
\]
and those of $d_\ell \Gamma\big(d(1-s)-e_\ell\big)$ with $0\leq \ell \leq h$ and $d_\ell\neq0$, at 
\[
s= 1+\frac{\nu}{d} \qquad \text{with} \  \nu\geq -e_\ell,
\]
must cancel. Hence for negative $\nu$ (if any) and $0\leq \ell\leq h$ we must have
\[
\overline{F}_\ell \big(\frac{\nu}{d},\alpha\big)=0 \qquad \text{for} \ \nu=-e_\ell,\dots,-1 \ \text{ whenever $d_\ell\neq0$},
\]
while for $\nu\geq0$ the residues have to cancel. But
\[
\begin{split}
\res_{s=1+\nu/d} \Gamma\big(d(1-s)-e_\ell\big) &= - \frac{1}{d} \frac{(-1)^{\nu+e_\ell}}{(\nu+e_\ell)!} \\
\res_{s=1+\nu/d} \Gamma\big(d(1-s)\big) &= - \frac{1}{d} \frac{(-1)^{\nu}}{\nu!},
\end{split}
\]
hence taking residues at $s=1+\nu/d$ with $\nu\geq0$ in \eqref{1-30}, and recalling \eqref{1-19} with $m_F=1$, we obtain
\begin{equation}
\label{4-20}
 \frac{\omega_F}{\sqrt{2\pi}} \left(\frac{q^{1/d}}{2\pi d}\right)^{-\frac{d}{2}-\nu} \sum_{\ell=0}^h \frac{d_\ell}{d} \frac{(-1)^{e_\ell}}{(\nu+e_\ell)!}  \overline{F}_\ell \big(- \frac{\nu}{d},\alpha\big) + (2\pi i \alpha)^\nu \rho_F = 0.
\end{equation}
Note that \eqref{4-20} holds for $\alpha\in$ Spec$(F)$ as well, since $F(s,\alpha)$ has no poles at the points $s=1+\nu/d$, $\nu\geq0$. In particular, \eqref{4-20} gives a formula for the residue $\rho_F$, of course under the above hypotheses.

\smallskip
{\bf 6.} There are some differences between the present approach to the functional equation of $F(s,\alpha)$ and that in \cite{Ka-Pe/half}. Indeed, in \cite{Ka-Pe/half} we were led directly to the strict functional equation, while the present method gives always \eqref{1-26} as a first instance, even if $F(s,\alpha)$ satisfies a strict functional equation. However, later on we somehow recover the previous approach, since the manipulations of the $\Gamma$ function in Section 4.2 correspond, in a general situation, to those used in \cite{Ka-Pe/half} to show that actually \eqref{3-12} and \eqref{3-13} are algebraic identities, without error terms.

\smallskip
{\bf 7.} Finally, recalling Remark 1.2, we give a couple of simple applications of the algorithm for the computation of the value of $h$ in the strict functional equation \eqref{1-30}. We first deal with \eqref{1-4}. In this case the $\gamma$-factor in (iii) of Theorem 4 has $N=1$, $d=2$ and $n_1=h+2$, where $h\geq0$. Then $2N+1-n_1 = 1-h$ and we have
\[
\begin{split}
1-h &= (1-h) + 0 +0, \quad \text{hence} \ \nu_0 = 0 \\
h+2 &= (1-h) + 1 +2h, \quad \text{hence} \ \nu_1=h.
\end{split}
\]
Thus $N_0=h$, hence the value of $h$ in \eqref{1-4} is confirmed by means of Theorem 4. Similarly, in the case of \eqref{1-5} we have $N=1$, degree $d$ and $n_1= h+2$ with $h\geq0$ if the sign $+$ occurs, and $n_1=h+1$ with $h\geq1$ if $-$ occurs. The first case is exactly as the previous one, while in the second we have $2N+1-n_1=2-h$, thus
\[
\begin{split}
h+1 &= (1-h) + 0 +2h, \quad \text{hence} \ \nu_0 = h \\
2-h &= (1-h) + 1 +0, \quad \text{hence} \ \nu_1=0.
\end{split}
\]
Therefore in both cases $N_0=h$, thus the strict functional equation holds with the value $h$ in these cases as well.

\bigskip
\section{Appendix}

\smallskip
Here we gather the proofs of the auxiliary results used in the previous sections. The first one is an explicit bound for the error term in the truncated binomial series; here $\HH$ denotes the upper half-plane. Most probably it is known in the literature, but we couldn't locate it.

\medskip
{\bf Lemma A.} {\sl Let $z\in\overline{\HH}$, $\rho\in\CC$, $0<\delta<1/2$ and $|z|<1-\delta$. Then for every integer $R\geq0$ we have
\[
(1+z)^\rho = \sum_{r=0}^R {\rho \choose r} z^r + Q_R(z,\rho)
\]
where the branch of $(1+z)^\rho$ is as in \eqref{1-21} and}
\begin{equation}
\label{6-1}
Q_R(z,\rho)   \ll_{R,\delta} |z \rho|^{R+1} \delta^{-|\Re\rho|}  \max_{0\leq u \leq 1} e^{-\Im\rho \arg(1+uz)}.
\end{equation}

\medskip
{\it Proof.} We use Taylor's formula with the integral form of the remainder, see e.g. Section 3.4 of Duren \cite{Dur/2012}. Writing $z=xe^{i\theta}$ we have
\[
\begin{split}
Q_R(z,\rho) &= \frac{1}{R!} \int_0^x \frac{\partial^{R+1}}{\partial t^{R+1}} (1+te^{i\theta})^\rho (x-t)^R \d t \\
&= e^{i(R+1)\theta} (R+1) {\rho\choose R+1} \int_0^x (1+e^{i\theta}t)^{\rho-R-1} (x-t)^R \d t \\
&=  z^{R+1} (R+1) {\rho\choose R+1}  \int_0^1 (1+uz)^{\rho-R-1} (1-u)^R \d u \\
&\ll_R |z \rho|^{R+1} \max_{0\leq u \leq 1} |1+uz|^{\Re\rho-R-1} e^{-\Im\rho \arg(1+uz)} \\
&\ll_{R,\delta}  |z \rho|^{R+1} \delta^{-|\Re\rho|}  \max_{0\leq u \leq 1} e^{-\Im\rho \arg(1+uz)},
\end{split}
\]
as claimed. \qed

\medskip
Let, for $A,B\geq0$, $A+B>1$ and $t\in\RR$,
\[
I= I(t,A,B) := \int_{-\infty}^{+\infty} \frac{\d v}{(1+|v+t|)^A(1+|v|)^B}.
\]

\medskip
{\bf Lemma B.} {\sl We have}
\[
I \ll \frac{\log(2+|t|)}{(2+|t|)^{\min(A,B,A+B-1)}}.
\]

\medskip
{\it Proof.} We may assume that $|t|\geq 2$, otherwise the assertion is trivial since the integral is clearly convergent. We split the integral as
\[
I =\left(  \int_{|v+t|\leq |t|/2}  + \int_{\substack{|v+t|>|t|/2 \\ |v|\leq 2|t|}} + \int_{\substack{|v+t|>|t|/2 \\ |v| > 2|t|}} \right) \frac{\d v}{(1+|v+t|)^A(1+|v|)^B} = I_1+I_2+I_3,
\]
say. In the range of $I_1$ we have $|v| = |-t+v+t|  \geq |t|-|v+t| \geq |t|/2$, hence
\[
\begin{split}
I_1 \ll \frac{1}{(1+|t|)^B} \int_{-|t|}^{|t|} &\frac{\d v}{(1+|v|)^A} \ll 
\begin{cases}
(2+|t|)^{-B} \log(2+|t|) & \ \text{if} \ A\geq 1 \\
(2+|t|)^{-A-B+1} \log(2+|t|) & \ \text{if} \ 0\leq A\leq 1 
\end{cases} \\
&\ll \frac{\log(2+|t|)}{(2+|t|)^{\min(B,A+B-1)}}.
\end{split}
\]
Similarly,
\[
I_2 \ll \frac{1}{(1+|t|)^A} \int_{-2|t|}^{2|t|} \frac{\d v}{(1+|v|)^B} \ll \frac{\log(2+|t|)}{(2+|t|)^{\min(A,A+B-1)}}.
\]
Finally, in the range of $I_3$ we have $|v+t| \geq |v|-|t| \geq |v|/2$, thus
\[
I_3 \ll \int_{2|t|}^{+\infty} \frac{\d v}{v^{A+B}} \ll \frac{1}{(2+|t|)^{A+B-1}},
\]
and the lemma follows. \qed

\medskip
{\bf Remark 6.1.} The above proof shows that if $A,B\neq1$, then the factor $\log(2+|t|)$ can be avoided in Lemma B. \qed

\medskip
Let now 
\begin{equation}
\label{6-2}
P(z) = 1 + \sum_{\ell=1}^h \frac{b_\ell}{(z-1)\cdots(z-\ell)},
\end{equation}
where $h\geq0$ and $b_\ell\in\CC$ (the sum equals 0 if $h=0$). Note that for
\[
|z|> R_0:= h+\sum_{\ell=1}^h |b_\ell|
\]
we have
\[
\Big| \sum_{\ell=1}^h \frac{b_\ell}{(z-1)\cdots(z-\ell)}\Big| \leq  \sum_{\ell=1}^h \frac{|b_\ell|}{(|z|-1)\cdots(|z|-\ell)} \leq \frac{1}{|z|-h} \sum_{\ell=1}^h |b_\ell| <1,
\]
hence both $P(z)$ and $\log P(z)$ are holomorphic for $|z|>R_0$. Moreover, since $P(z)\to1$ as $|z|\to\infty$, for $|z|>R_0$ we have that
\begin{equation}
\label{6-3}
\log P(z)=\sum_{m=1}^\infty \frac{c_m}{z^m}.
\end{equation}
Thus, by Cauchy's formula applied to $\log P(1/z)$, we see that for $0<r<1/R_0$ and $m\geq1$
\[
|c_m| \leq \frac{1}{r^m} \max_{|z|=r} |\log P(1/z)|,
\]
therefore
\[
c_m \ll R_1^m
\]
for some fixed $R_1>R_0$. Consider now the entire function
\begin{equation}
\label{6-4}
\psi_P(w) = \sum_{m=1}^\infty \frac{c_m}{(m-1)!} w^m,
\end{equation}
satisfying
\begin{equation}
\label{6-5}
\psi_P(w) \ll \sum_{m=1}^\infty \frac{|R_1w|^m}{(m-1)!} \ll |w| e^{R_1|w|}.
\end{equation}

\medskip
{\bf Lemma C.} {\sl Let $P(z)$ and $\psi_P(w)$ be as in \eqref{6-2} and \eqref{6-4}, respectively, $R_1$ be as in \eqref{6-5} and let $0<\delta<\pi$ and $\eta>0$. Then for $|z|\geq cR_1$, with a certain $c=c(\delta)>0$, and $|\arg(z)| <\pi-\delta$ we have
\[
\log P(z) = \int_0^{\eta e^{i\beta}} \psi_P(w) e^{-zw} \frac{\d w}{w} + O_{\delta,\eta}\big(e^{-\kappa|z|}\big)
\]
with a certain $\kappa=\kappa(\delta,\eta)>0$, where $|\beta|<(\pi-\delta)/2$ is such that $|\arg(z)+\beta| < (\pi-\delta)/2$.}

\medskip
{\it Proof.} Write
\[
I(z) = \int_0^{\infty e^{i\beta}} \psi_P(w) e^{-zw} \frac{\d w}{w}.
\]
In view of \eqref{6-5} and the hypotheses of the lemma, for $|w|\geq \eta$ the integrand is
\begin{equation}
\label{6-6}
\ll e^{R_1|w| - \Re(zw)} \ll e^{R_1|w| - |zw|\cos((\pi-\delta)/2)} \ll e^{-\kappa_0 |wz|}
\end{equation}
with a certain $\kappa_0=\kappa_0(\delta)>0$, provided $|z|\geq cR_1$ with a certain $c=c(\delta)$. Moreover, as $w\to0$ we have $\psi_P(w) = O(|w|)$, hence $I(z)$ is absolutely convergent provided
\begin{equation}
\label{6-7}
|\arg(z)+\beta| < (\pi-\delta)/2 \quad \text{and} \quad |z|\geq cR_1.
\end{equation}

\smallskip
Next we compute $I(z)$. From \eqref{6-4}, condition $|\arg(z)+\beta| < (\pi-\delta)/2$ and Cauchy's theorem, thanks to the decay of the exponential function we have
\[
\begin{split}
I(z) = \sum_{m=1}^\infty \frac{c_m}{(m-1)!} \int_0^{\infty e^{i\beta}} &e^{-zw} w^{m-1} \d w =  \sum_{m=1}^\infty \frac{c_m}{z^m(m-1)!} \int_0^{\infty e^{i(\beta + \arg(z))}} e^{-w} w^{m-1} \d w \\
&=  \sum_{m=1}^\infty \frac{c_m \Gamma(m)}{z^m(m-1)!} = \log P(z)
\end{split}
\]
by \eqref{6-3}. Moreover, from \eqref{6-6} we get
\[
 \int_{\eta e^{i\beta}}^{\infty e^{i\beta}} \psi_P(w) e^{-zw} \frac{\d w}{w} \ll  \int_{\eta e^{i\beta}}^{\infty e^{i\beta}} e^{-\kappa_0 |wz|} |\d w| \ll_{\delta,\eta} e^{-\kappa_0\eta |z|}
\]
provided \eqref{6-7} holds, and the lemma follows. \fine

\medskip
Finally, we need the following lemma, which (essentially) was stated without proof as Lemma 2.2 in \cite{Ka-Pe/2002a}. For $a\in\CC$ we write
\[
\psi(w,a) = \sum_{n=2}^\infty \frac{(-1)^n B_n(a)}{n!} w^{n-1}, \hskip1.5cm |w|<2\pi.
\]

\medskip
{\bf Lemma D.} {\sl Let $\gamma,\delta,\beta\in\RR$ satisfy
\[
0<\gamma\leq 1, \quad 0<\delta<\half \quad \text{and} \quad |\beta|<\pi(\half-\delta).
\]
Then for $\rho,a\in\CC$ satisfying
\[
|\rho|\geq \frac{1}{\delta} (1+|a|) \quad \text{and} \quad |\arg(\rho)+\beta| <\pi(\half-\delta)
\]
we have
\[
\begin{split}
\log\Gamma(\rho+a) =& (\rho+a-\half)\log\rho -\rho +\half\log(2\pi) \\
&+ \int_0^{\gamma e^{i\beta}} \psi(w,a) e^{-\rho w} \frac{\d w}{w} +O_\gamma\big(\frac{1}{\delta} e^{-\gamma\delta|\rho|/2}\big)
\end{split}
\]
uniformly in $a$ and $\beta$.}

\medskip
{\it Proof.} We start with the expression
\begin{equation}
\label{6-8}
\log \Gamma(\rho) = (\rho-\half) \log\rho - \rho + \half\log(2\pi) +\int_0^{\infty e^{i\beta}} \psi(w) e^{-w\rho} \frac{\d w}{w},
\end{equation}
where
\[
\psi(w) = \frac{1}{e^w-1} - \frac{1}{w} +\frac{1}{2}
\]
and the formula is valid for $|\beta|<\pi/2$, $|\arg(\rho)|<\pi$ and $|\arg(\rho)+\beta|<\pi/2$; see (5) of Section 1.9 of \cite{EMOT/1953-I}. We apply \eqref{6-8} with $\rho$ replaced by $\rho+a$, hence we have to check the conditions in this case. From our hypotheses on $\rho,a,\beta$ and $\delta$, by elementary geometry we have
\[
\begin{split}
|\arg(\rho+a)| &= |\arg(\rho) + \arg(1+a/\rho)| \leq |\arg(\rho)+\beta| +|\beta| + \arg(1+a/\rho)| \\
&\leq \pi(\half-\delta) + \pi(\half-\delta) + \arcsin|a/\rho| \\
&\leq \pi(1-2\delta) + \frac{\pi}{2}|a/\rho| \leq \pi(1-\frac{3}{2}\delta)
\end{split}
\]
and, similarly,
\begin{equation}
\label{6-9}
|\arg(\rho+a)+\beta| \leq \pi(\half-\delta) + \frac{\pi}{2}|a/\rho| \leq \frac{\pi}{2} (1-\delta).
\end{equation}
Thus from \eqref{6-8} we get
\begin{equation}
\label{6-10}
\begin{split}
\log \Gamma(\rho+a) &= (\rho+a-\half) \log\rho - \rho + \half\log(2\pi) \\
&+ \ell(\rho,a) +\int_0^{\infty e^{i\beta}} \psi(w) e^{-w(\rho+a)} \frac{\d w}{w},
\end{split}
\end{equation}
where
\begin{equation}
\label{6-11}
\ell(\rho,a) = (\rho+a-\half)\log(1+a/\rho) -a.
\end{equation}

\smallskip
Since $w=re^{i\beta}$ with $r>0$, thanks to \eqref{6-9} we have
\begin{equation}
\label{6-12}
\begin{split}
|e^{-w(\rho+a)}| &= e^{-r|\rho+a| \cos(\arg(w(\rho+a)))} \leq e^{-r|\rho|(1-|a/\rho|) \cos(\frac{\pi}{2}(1-\delta))} \\
&\ll e^{-r|\rho|(1-\delta) \sin(\pi\delta/2)} \ll e^{-r|\rho|\delta/2},
\end{split}
\end{equation}
and for $r\geq \gamma$
\[
|\psi(w)| \leq \frac{1}{|e^w-1|} + \frac{1}{|w|} + \frac12 \ll_\gamma \frac{1}{\delta}.
\]
Hence
\begin{equation}
\label{6-13}
\int_{\gamma e^{i\beta}}^{\infty e^{i\beta}} \psi(w) e^{-w(\rho+a)} \frac{\d w}{w} \ll_\gamma \frac{1}{\delta} \int_\gamma^\infty e^{-r|\rho|\delta/2} \frac{\d r}{r} \ll_\gamma \frac{1}{\delta^2|\rho|} e^{-\gamma\delta|\rho|/2} \ll_\gamma \frac{1}{\delta} e^{-\gamma\delta|\rho|/2}.
\end{equation}

\smallskip
Next we have
\begin{equation}
\label{6-14}
\begin{split}
\psi(w,a) &= \frac{1}{w} \Big(\sum_{n=0}^\infty \frac{B_n(a)}{n!} (-w)^n -1+B_1(a)w\big) \\
&=  \frac{1}{w} \Big(\frac{-we^{-aw}}{e^{-w}-1} -1 +(a-\half)w\Big) = \frac{e^{(1-a)w}}{e^w-1} -\frac{1}{w} +a-\frac12 \\
&= \frac{e^{-aw}(e^w-1+1)}{e^w-1}  -\frac{1}{w} +a-\half = \frac{e^{-aw}}{e^w-1} +e^{-aw}  -\frac{1}{w} +a-\frac12 \\
&= e^{-aw}\big(\psi(w) +\frac{1}{w} -\frac12\big) +e^{-aw}  -\frac{1}{w} +a-\frac12 = e^{-aw}\psi(w) + h(w,a),
\end{split}
\end{equation}
where
\[
h(w,a) = \frac{e^{-aw}}{w} + \frac{e^{-aw}}{2} -\frac{1}{w} +a-\frac12.
\]
Note that
\begin{equation}
\label{6-15}
h(w,a) \ll_a |w| \quad \text{as} \ \ w\to0
\end{equation}
since
\[
h(w,a) = \frac{1}{w} \big(1-aw+O_a(|w|^2)\big) + \frac12 \big(1+O_a(|w|)\big) -\frac{1}{w} +a-\frac12 \ll_a |w|.
\]
Moreover, for $w=re^{i\beta}$ with $r\geq \gamma$, by \eqref{6-12} and our hypotheses we have
\begin{equation}
\label{6-16}
h(w,a)e^{-\rho w} \ll_\gamma |e^{-w(\rho+a)}| +(1+|a|)|e^{-\rho w}| \ll_\gamma (1+|a|) e^{-r|\rho|\delta/2}.
\end{equation}
Hence
\begin{equation}
\label{6-17}
\begin{split}
\int_{\gamma e^{i\beta}}^{\infty e^{i\beta}} h(w,a) e^{-\rho w} \frac{\d w}{w} &\ll_\gamma (1+|a|) \int_\gamma^\infty e^{-r|\rho|\delta/2} \frac{\d r}{r} \\
&\ll_\gamma \frac{1+|a|}{\delta|\rho|} e^{-\gamma\delta|\rho|/2} \ll_\gamma e^{-\gamma\delta|\rho|/2}.
\end{split}
\end{equation}
Thus, from \eqref{6-14} and \eqref{6-17} we get
\begin{equation}
\label{6-18}
\int_0^{\gamma e^{i\beta}} \psi(w) e^{-(\rho+a)w} \frac{\d w}{w} = \int_0^{\gamma e^{i\beta}} \psi(w,a) e^{-\rho w} \frac{\d w}{w}  - H(\rho,a) + O_\gamma\big(e^{-\gamma\delta|\rho|/2}\big),
\end{equation}
where
\[
H(\rho,a) = \int_0^{\infty e^{i\beta}} h(w,a) e^{-\rho w} \frac{\d w}{w}.
\]
Note that, thanks to \eqref{6-15} and \eqref{6-16}, this integral is absolutely and uniformly convergent on compact sets of $\rho$ satisfying our hypotheses. To conclude the proof we now show that 
\begin{equation}
\label{6-19}
H(\rho,a) = \ell(\rho,a),
\end{equation}
where $\ell(\rho,a)$ is defined by \eqref{6-11}. Indeed, the lemma follows at once from \eqref{6-10},\eqref{6-13},\eqref{6-18} and \eqref{6-19}.

\smallskip
From the definition of $H(\rho,a)$ and $h(w,a)$ we have that
\[
\begin{split}
\frac{\partial^2}{\partial\rho^2} H(\rho,a) &= \int_0^{\infty e^{i\beta}} h(w,a) e^{-\rho w} w \d w = \int_0^{\infty e^{i\beta}} e^{-(\rho+a)w} \d w \\
&+ \frac12 \int_0^{\infty e^{i\beta}} w e^{-(\rho+s)w} \d w - \int_0^{\infty e^{i\beta}} e^{-\rho w}\d w + (a-\half) \int_0^{\infty e^{i\beta}} w e^{-\rho w} \d w \\
&= \frac{1}{\rho+a} + \frac12 \frac{1}{(\rho+a)^2} - \frac{1}{\rho} + (a-\half) \frac{1}{\rho^2},
\end{split}
\]
while trivially from the defnition of $\ell(\rho,a)$ we have
\[
\frac{\partial^2}{\partial\rho^2}  \ell(\rho,a) = \frac{1}{\rho+a} + \frac12 \frac{1}{(\rho+a)^2} - \frac{1}{\rho} + (a-\half) \frac{1}{\rho^2}.
\]
Thus
\[
H(\rho,a) = \ell(\rho,a) +A\rho +B
\]
for some constants $A$ and $B$. But from \eqref{6-11} we see that as $|\rho|\to\infty$
\[
\ell(\rho,a) \ll_a \frac{1}{|\rho|},
\]
while from our hypotheses, \eqref{6-15} and \eqref{6-17} we get
\[
\begin{split}
H(\rho,a) &\ll_\gamma \int_0^{\gamma e^{i\beta}} h(w,a) e^{-\rho w} \frac{\d w}{w} + O_\gamma \big(e^{-\gamma\delta|\rho|/2}\big) \\
&\ll_{\gamma,a} \int_0^\gamma e^{-\delta|\rho|r} \d r + e^{-\gamma\delta|\rho|/2}. \ll_{\gamma,a} \frac{1}{\delta|\rho|}  e^{-\gamma\delta|\rho|/2} \ll_{\gamma,a,\delta} \frac{1}{|\rho|}.
\end{split}
\]
Therefore
\[
A\rho+B \ll \frac{1}{|\rho|}
\]  
and hence $A=B=0$. Thus \eqref{6-19} holds and the lemma is proved. \qed

\bigskip

\ifx\undefined\bysame{poly}.
\newcommand{\bysame}{\leavevmode\hbox to3em{\hrulefill}\ ,}
\fi

\bigskip
\bigskip
\bigskip
\noindent
Jerzy Kaczorowski, Faculty of Mathematics and Computer Science, A.Mickiewicz University, 61-614 Pozna\'n, Poland and Institute of Mathematics of the Polish Academy of Sciences, 
00-956 Warsaw, Poland. e-mail: kjerzy@amu.edu.pl

\medskip
\noindent
Alberto Perelli, Dipartimento di Matematica, Universit\`a di Genova, via Dodecaneso 35, 16146 Genova, Italy. e-mail: perelli@dima.unige.it

\end{document}